\documentclass[a4paper,leqno,11pt,twoside]{amsart}
\usepackage{amsmath,amsfonts,amssymb,amsthm,amscd}
\usepackage[utf8]{inputenc}
\usepackage[english]{babel}
\usepackage[colorlinks=true,citecolor=blue, urlcolor=blue, linkcolor=blue,pagebackref]{hyperref}
\usepackage[top=1.1in,bottom=1.1in,left=1.in,right=1.in]{geometry} 
\usepackage{graphicx}
\usepackage{wrapfig}
\usepackage{paralist}
\usepackage{tabto}
\usepackage{standalone}
\usepackage{xfrac}
\usepackage{float}
\usepackage{tikz-cd}
\usepackage{mathtools}

\usepackage{pgf,tikz}
\usetikzlibrary{arrows,chains,
 decorations.pathreplacing,
shapes,%
}

\usepackage[all]{xy} 

\usepackage{enumerate}
\usepackage{xcolor}
\usepackage{aliascnt}

\theoremstyle{plain}
\newtheorem{thm}{Theorem}[section]

\newtheorem{thmIntr}{Theorem}

\newaliascnt{propIntr}{thmIntr}

\newaliascnt{corIntr}{thmIntr}
\newtheorem{corIntr}[corIntr]{Corollary}

\newaliascnt{QU}{thm}
\newtheorem{QU}[QU]{Question}
\aliascntresetthe{QU}

\newaliascnt{lem}{thm}
\newtheorem{lem}[lem]{Lemma}
\aliascntresetthe{lem}

\newaliascnt{cor}{thm}
\newtheorem{cor}[cor]{Corollary}
\aliascntresetthe{cor}

\newaliascnt{prop}{thm}
\newtheorem{prop}[prop]{Proposition}
\aliascntresetthe{prop}

\theoremstyle{definition}

\newaliascnt{rem}{thm}
\newtheorem{rem}[rem]{Remark}
\aliascntresetthe{rem}

\newaliascnt{defn}{thm}

\aliascntresetthe{defn}

\newaliascnt{ex}{thm}

\aliascntresetthe{ex}

\numberwithin{equation}{section}

\def\bP{\ensuremath{\mathbb{P}}}
\def\bQ{\ensuremath{\mathbb{Q}}}
\def\bR{\ensuremath{\mathbb{R}}}
\def\bZ{\ensuremath{\mathbb{Z}}}
\def\bC{\ensuremath{\mathbb{C}}}

\def\cE{\ensuremath{\mathcal{E}}}
\def\cF{\ensuremath{\mathcal{F}}}
\def\cG{\ensuremath{\mathcal{G}}}
\def\cH{\ensuremath{\mathcal{H}}}
\def\cI{\ensuremath{\mathcal{I}}}
\def\cM{\ensuremath{\mathcal{M}}}
\def\cO{\ensuremath{\mathcal{O}}}

\def\cT{\ensuremath{\mathcal{T}}}

\def\Db{\mathop{\mathrm{D}^{\mathrm{b}}}\nolimits}
\DeclareMathOperator{\td}{td}
\DeclareMathOperator{\Pic}{Pic}
\DeclareMathOperator{\Supp}{supp}
\DeclareMathOperator{\Ext}{Ext}
\DeclareMathOperator{\ext}{ext}
\DeclareMathOperator{\sExt}{\mathcal{E}\emph{xt}}
\DeclareMathOperator{\Hom}{Hom}
\DeclareMathOperator{\sHom}{\mathcal{H}\emph{om}}
\DeclareMathOperator{\NS}{NS}
\DeclareMathOperator{\rk}{rk}

\DeclareMathOperator{\im}{Im}
\DeclareMathOperator{\ch}{ch}
\DeclareMathOperator{\Coh}{Coh}

\DeclareMathOperator{\coker}{coker}

\DeclareMathOperator{\irr}{irr}
\DeclareMathOperator{\Gr}{Gr}

\DeclareMathOperator{\colength}{colength}
\DeclareMathOperator{\length}{length}

\newcommand{\commr}[1]{{\color{red}#1}} 
\definecolor{applegreen}{rgb}{0.55, 0.71, 0.0}
\newcommand{\commg}[1]{{\color{applegreen}#1}} 

\newcommand{\set}[1]{\left\{#1\right\}}

\title[On the degree of irrationality of low genus $K3$ surfaces]{On the degree of irrationality of low genus $K3$ surfaces}
\author[F.~Moretti and A.~Rojas]{Federico Moretti and Andr\'es Rojas}

\address{F.M.: Institut für Mathematik, Humboldt-Universität zu Berlin, Unter den Linden 6, 10099 Berlin, Germany}
\curraddr{Stony Brook University,  100 Nicolls Road, Stony Brook, NY 11794}
\email{federico.moretti@stonybrook.edu} 

\address{A.R.: Institut für Mathematik, Humboldt-Universität zu Berlin, Unter den Linden 6, 10099 Berlin, Germany}
\email{rojasand@math.hu-berlin.de} 

\begin{document}

\begin{abstract}
Given a general polarized $K3$ surface $S\subset \mathbb P^g$ of genus $g\le 14$, we study projections $S\hookrightarrow \mathbb P^g\dasharrow \mathbb P^2$ of minimal degree and their variational structure. In particular, we prove that the degree of irrationality of all such surfaces is at most $4$, and that for $g=7,8,9,11$ there are no rational maps $S\dasharrow\bP^2$ of degree $3$ induced by the primitive linear system. Our methods combine vector bundle techniques à la Lazarsfeld with derived category tools, and also make use of the rich theory of singular curves on $K3$ surfaces.
\end{abstract}


\maketitle

\setcounter{tocdepth}{1}

\section{Introduction}

Given an irreducible, projective variety $X$ of dimension $n$, its \emph{degree of irrationality} $\irr(X)$ is the minimal degree of a rational dominant map $X\dasharrow \bP^n$. This is a birational invariant measuring how far $X$ is from being rational, which in recent times has been the object of a considerable amount of work; remarkable progress has been made for hypersurfaces in projective space (\cite{BDELU,fanohyp}), abelian varieties (\cite{nathchen,martin,av}), hyperkähler varieties and their moduli spaces (\cite{hk,moduliK3,moduliHK}), and some concrete examples (\cite{c3f}). Nevertheless, the computation of $\irr(X)$ is still a difficult problem for which very little is known in general if $n\geq 2$.

A preliminary version of this problem consists of fixing a line bundle $L\in\Pic(X)$, and considering the following invariant of the pair $(X,L)$, closely related to the degree of irrationality of $X$:
\[
\irr_L(X):=\min\{\deg\left(\varphi_V:X\dasharrow\bP^n\right)|\varphi_V\text{  dominant defined by $V\in\Gr(n+1, H^0(L))$}\}
\]

The degree of a dominant map $\varphi_V:X\dasharrow \bP(V^\vee) = \bP^n$ has a natural interpretation in terms of the kernel $E^\vee$ of the evaluation map $V \otimes \cO_X\to L$. Indeed, the inclusion $V^\vee\subset H^0(E)$ describes the fibers of $\varphi_V$ in a precise way (see \autoref{compdegree} for the surface case). The degree of $\varphi_V$ turns out to be the difference between $c_n(E)$ and the degree of the ``cycle-theoretic base locus" of sections of $V^\vee\subset H^0(E)$. Note that the rational map induced by $V\in\Gr(n+1, H^0(L))$ is nothing but the composition
\[
X\dasharrow\bP(H^0(L)^\vee)\dasharrow \bP(V^\vee),
\]
where the second map is the projection from the linear space $\bP\left(\ker(H^0(L)^\vee\to V^\vee)\right)$. In the case of surfaces, the top Chern class $c_2(E)$ gives information on \emph{how secant} is this linear space with respect to $X$, whereas the ``cycle-theoretic base locus" says \emph{how tangent} is this linear space. If one is able to control these two invariants efficiently, then one can infer on the projections of a given degree.

This approach has been exploited by the first author in \cite{mor}, in order to obtain new upper bounds for the quantities $\irr_L(X)$ and $\irr(X)$ for several examples of polarized varieties $(X,L)$. 
Among those examples, one finds polarized $K3$ surfaces $(S,L)$ with $\Pic(S)=\bZ\cdot L$ (see \cite[Theorem 1.1]{mor}). In that case the kernel bundle is moreover stable, which allowed the first author  to endow the following set of special maps with a natural scheme structure:
\[
W^2_d(S,L):=\{V\in\Gr(3, H^0(L))\;|\;\varphi_V:S\dasharrow \bP^2\text{ is gen. finite of degree $\leq d$}\}
\]

The goal of this paper is to show how, in the case of $K3$ surfaces, a combination of vector bundle techniques with tools coming from the derived category (mainly Bridgeland stability and Fourier-Mukai theory) provides new insight into these problems.

More precisely, we examine in detail the invariants $\irr_L(S)$ and $W^2_{\irr_L(S)}(S,L)$ for polarized $K3$ surfaces $(S,L)$ of low genus $g$, under the assumption $\Pic(S)=\bZ\cdot L$. This question is already known for $g\leq 6$ (see \cite[Theorem 1.2]{mor}); as it will become apparent in our arguments, the problem is substantially harder as long as $g$ increases, due to the complexity of the potential cycle-theoretic base loci.

Since the picture is different depending on the parity of $g$, we will state our results separately. For a polarized $K3$ surface $(S,L)$ of odd genus $g$, let us denote by $\cM$ the moduli space of rank $2$ stable vector bundles with $c_1=L$ and minimal $c_2=\lfloor\frac{g+3}{2}\rfloor$. 

\begin{thmIntr}\label{thmA}
    Let $(S,L)$ be a very general polarized $K3$ surface of genus $g$.
\begin{enumerate}[{\rm (1)}]
        \item If $g=7$, then $\irr_L(S)=4$ and $W^2_4(S,L)$  is isomorphic to $S\times \cM$.
        \item If $g=9$, then $\irr_L(S)=4$ and $W^2_4(S,L)$ has at least two irreducible components: a $3$-dimensional one (isomorphic to a $\mathbb P^1$-bundle over $\cM$) and a $2$-dimensional one (isomorphic to a nontrivial correspondence between $S$ and $\cM$). 
        \item If $g=11$, then $\irr_L(S)=4$ and $W^2_4(S,L)$ has an irreducible component isomorphic to $\cM$.
        \item If $g=13$, then $\irr_L(S)\leq4$ and $W^2_4(S,L)$ is at least $1$-dimensional.
    \end{enumerate}
    \label{A}
\end{thmIntr}

For even genera, our results can be summarized as follows:

\begin{thmIntr}\label{thmB}
    Let $(S,L)$ be a polarized $K3$ surface of genus $g$, with $\Pic(S)=\bZ\cdot L$.
    \begin{enumerate}[{\rm (1)}]
        \item If $g=8$, then $\irr_L(S)=4$ and $W^2_4(S,L)$  has a component birational  to a $\mathbb P^3$-bundle over $S$.
        \item If $g=10$, then $\irr_L(S)\le 4$ and $W^2_4(S,L)$ has at least two irreducible components: a $4$-dimensional one (isomorphic to the Hilbert square $S^{[2]}$) and a $3$-dimensional one (isomorphic to $\mathbb P(E)$, where $E$ is the stable spherical rank $2$ vector bundle with $c_2(E)=6$).
        \item If $g=12$, then $\irr_L(S)\le 4$ and $W^2_4(S,L)$ has a unirational $3$-dimensional component.
        \item If $g=14$, then $\irr_L(S)\le 4$ and $W^2_4(S,L)$ has an irreducible component isomorphic to $S$.
    \end{enumerate}
    \label{B}
\end{thmIntr}

Since upper bounds for the degree of irrationality specialize in families of regular surfaces (see \cite[Proposition~C]{fanohyp}), an immediate consequence of \autoref{thmA} and \autoref{thmB} is:

\begin{corIntr}
Any polarized $K3$ surface of genus $\leq 14$ has degree of irrationality $\leq 4$.
\end{corIntr}

Let us say a few words about the strategy of proof of \autoref{thmA}, which is the main body of work. The first step is the construction of rational maps $S\hookrightarrow\bP(H^0(L)^\vee)\dasharrow\bP^2$ of degree $\leq 4$; by \autoref{compdegree}, this can be achieved by finding couples $(E,\xi)\in\cM\times S^{[\frac{g-5}{2}]}$ such that $h^0(E\otimes\cI_\xi)\geq 3$. Via the derived category approach, we determine the locus of such couples, which leads to an \emph{explicit} description of the corresponding irreducible component of $W^2_4(S,L)$.

For instance, in genus 11 Bridgeland stability provides an isomorphism of Hilbert schemes
\[
\gamma:\cM^{[3]}\overset{\cong}{\longrightarrow} S^{[3]}
\]
that, combined with Fourier-Mukai arguments, shows that
\[
\{(E,\xi)\in\cM\times S^{[3]}\;|\;h^0(E\otimes\cI_\xi)\geq 3\}=\{(E,\gamma(\eta_E))\;|\;E\in\cM\}\cong\cM
\]
(here $\eta_E$ denotes the non-curvilinear element of $\cM^{[3]}$ supported at $E$).

In the second step (for $g=7,9,11$) we prove that there are no degree 3 maps; this is divided into two substeps. One starts by discarding the existence of couples $(E,\xi)\in\cM\times S^{[ \frac{g-3}{2}]}$ such that $h^0(E\otimes\cI_\xi)\geq 3$, again via derived methods. Geometrically, this means discarding the situation where the linear space we are projecting from is maximally secant to $S$ (i.e. the base locus has the highest possible length) and is tangent to $S$ at $\frac{g-3}{2}$ points. For instance in genus $g=9$, by (twisted) Fourier-Mukai theory the condition $h^0(E\otimes \cI_\xi)\geq 3$ implies the existence of a rational curve in $\cM$ with a triple point; coupled with a result of Chen \cite{chen} on the singularities of rational curves on general $K3$ surfaces, this gives a contradiction.

Then one has to rule out all the cases where the ``cycle-theoretic base locus" is strictly bigger than the ``scheme-theoretic base locus". This is a more special situation, in which the linear space we are projecting from is maximally secant and the singularities of the intersection with $S$ are more complicated. In this case, we exploit systematically an \emph{ad hoc} lemma relating the ``cycle-theoretic base locus" with the image of $V^\vee\otimes \cO_S\to E$ (c.f. \autoref{ideallemma}), which we combine with other vector bundle arguments.

It is worth mentioning that our results suggest several patterns for the behaviour of $\irr_L(S)$, as $g$ varies; here the efforts to describe $W^2_{\irr_L(S)}(S,L)$ become crucial. We thus expect these two invariants to be a potential source of new questions. Some of them are gathered in a separate section at the end of the paper.

Nevertheless, it seems to us that a good understanding of this problem for arbitrary $g$ is not accessible without further improvement of our techniques. As mentioned above, the reason is that the possible configurations of certain $0$-dimensional subschemes become more and more intricate. The case of genus $14$ is somehow the first concrete realization of this problem: maps of degree $4$ necessarily have a ``cycle-theoretic base locus" bigger than the ``scheme-theoretic base locus". In more geometric words, there is no codimension $3$ linear subspace in $\mathbb P^{14}$ intersecting $S$ along a subscheme of colength $18$, and being tangent to $S$ at four distinct points.

\vspace{3mm}

\textbf{Structure of the paper.} We start with a section of preliminaries, where we recall the kernel bundle approach to study linear systems (including Brill-Noether loci and key lemmas involving the local structure of the base ideal), as well as some basics on the derived category tools. Sections~\ref{sec:genus7}-\ref{sec:genus13} deal with the cases of odd genus (namely with the proof of \autoref{thmA}) and constitute the core of the paper. In section~\ref{sec:evengen} we treat even genera (i.e. \autoref{thmB}). Finally, in section \ref{sec:questions} we raise several questions about the behaviour of $\irr_L(S)$ and $W^2_{\irr_L(S)}(S,L)$ in arbitrary genus, based on our results for $g\leq 14$.

\textbf{Acknowledgements.} We would like to thank G. Farkas for mentoring, as well as A. Beauville, D. Huybrechts, M. Lahoz, E. Macrì, R. Moschetti, I. Smirnov and P. Stellari for answering some questions of ours. We are also grateful to R. Lazarsfeld, Á. Ríos Ortiz and I. Vogt for some interesting conversations, and to an anonymous referee for numerous comments and suggestions.

The authors are supported by the ERC Advanced Grant SYZYGY. This project has received funding from the European Research Council
	(ERC) under the European Union Horizon 2020 research and innovation program (grant agreement No. 834172).

\section{Preliminaries}\label{sec:prelim}

Throughout this paper we will work over the field of the complex numbers.

\subsection{Projections of low degree}
We fix a smooth polarized surface $(S,L)$ with $\mathrm{NS}(S)=\mathbb Z\cdot [L]$.
We will study linear systems $V\in \mathrm{Gr}(3,H^0(L))$, via the associated kernel bundle. Namely one can consider the vector bundle $E^\vee$ sitting in an exact sequence
\[
\begin{tikzcd}
    0 \arrow{r} & E^\vee \arrow{r} & V \otimes \mathcal O_S \arrow{r}{ev} & L \otimes \mathcal I \arrow{r} & 0,
\end{tikzcd}
\]
and let $T:=\underline{\mathop{\mathrm{Spec}}}(\cO_S/\cI)$ denote the base locus of the rational map $\varphi_V:S \dashrightarrow \mathbb P(V^\vee)$. Dualizing one gets an inclusion $V^\vee \subset H^0(E)$ that allows to study the fibers of $\varphi_V$ (note that $\varphi_V$ is dominant, by the assumption $\NS(S)=\mathbb Z\cdot [L]$). More precisely, the following holds:

\begin{prop}[{\cite[Section~2]{mor}}]\label{compdegree}
The morphism $\tilde \varphi={\varphi_V}|_{S\setminus T}:S\setminus T\to \bP(V^\vee)$ satisfies:

\begin{enumerate}[{\rm (1)}]
    \item For every point $p=[s_p]\in\bP (V^\vee)$, one has $\widetilde{\varphi}^{-1}(p)=Z(s_p)\cap (S\setminus T)$ where $s_p\in V^\vee\subset H^0(E)$.

    \item $\deg(\varphi)=c_2(E)-\deg\left(\bigcap_{s\in V^\vee}Z_{cycle}(s)\right)$.
\end{enumerate}
\end{prop}

The first preliminary lemma tells us that the degree of the map cannot be too big with respect to the colength of the base ideal $\cI$:

\begin{lem}\label{prellemma}
 Let $(S,L)$ be a polarized surface and $V\in \mathrm{Gr}(3,H^0(L))$ a linear system with finite base locus $\underline{\mathop{\mathrm{Spec}}}(\cO_S/\cI)$. Then $\mathrm{deg}(\varphi_V)\ge L^2- \frac{4}{3}\mathrm{colength}(\cI)$.
\end{lem}
\begin{proof}
There is another standard way to compute the degree of the rational map $\varphi_V$, in terms of Hilbert-Samuel multiplicities of $\cI$. The Hilbert-Samuel multiplicity of an ideal sheaf $\cI$ at a point $p$ is defined as $e_p(\cI):=\mathrm{length}(\cO_{S,p}/(f,g))$, where $f,g\in \cI\otimes \cO_{S,p}$ are general elements. Then \[
 \mathrm{deg}(\varphi_V)=L^2-\sum_{p \in S}e_p(\cI).
 \]
 Now by \cite[Proposition 5.7]{hmqs} $e_p(\cI)\le \frac{4}{3}\mathrm{colength}(\cI,p)$ (notice that $\cI\otimes \cO_{S,p}$ has at most three generators). The lemma follows.
 \end{proof}
 
 This has the following important implication, if one wants to study maps of degree $3,4$ in the primitive linear system of a $K3$ surface via the kernel bundle.
 
 \begin{cor}\label{minimalc2}
If $\varphi_V:S \dashrightarrow \mathbb P ^2$ is a map of degree $d$ with kernel bundle $E^\vee$, then
	\[
	c_2(E)\le \frac 1 4 (3d+L^2). \]
	In particular, for $(S,L)$ a polarized $K3$ surface of genus $g$ with $\Pic(S)=\bZ\cdot L$: if $d=3$ only the minimal $c_2=\lfloor \frac {g+3}{2} \rfloor$ is possible, and if $d=4$ only the two lowest $c_2$ are possible.
\end{cor}
\begin{proof}
The inequality follows from the previous lemma, together with the observation that $c_2(E)=L^2-\mathrm{colength}( \mathcal I)$.
 For the second part, observe that under the assumption $\Pic(S)=\bZ\cdot L$ the kernel bundle $E^\vee$ is stable; otherwise, we would have $\cO_S\hookrightarrow E^\vee$ which contradicts the inclusion $V\hookrightarrow H^0(L\otimes\cI)$.
 Hence by Mukai's theory of stable vector bundles on $K3$ surfaces (\cite[Corollary 2.5]{mvbnd}) we have $v(E)^2\geq -2$ (see \autoref{DerCatSec} for the definition of $v(E)$), which in this case reads as $c_2(E)\geq \lfloor \frac {g+3}{2} \rfloor$.
\end{proof}

In other words: if $S$ is a $K3$ surface with $\Pic(S)=\bZ\cdot L$ and $\varphi_V$ is of degree 3, then the base locus $T$ has the highest possible length (i.e. the codimension 3 linear subspace $\bP\left(\ker(H^0(L)^\vee\to V^\vee)\right)$ is \emph{maximally secant} to $S$).

\subsection{Brill-Noether loci}
Let $(S,L)$ be a polarized $K3$ surface of genus $g$ (i.e. $L^2=2g-2$), with $\Pic(S)=\bZ\cdot L$. In this subsection we recall the construction of the Brill-Noether loci from \cite{mor}, and prove some auxiliary lemmas that will be fundamental to study them from our perspective.

The set-theoretic definition is the following:
\[
W^2_d(S,L):=\{V \in \mathrm{Gr}(3,H^0(L))\; | \;\mathrm{deg}(\varphi_V)\le d\}.
\]
Let us denote by $\cM(2,L,g+1-c)$ the moduli space of stable rank 2 vector bundles with Chern classes $c_1=L$ and $c_2=c$ (this notation is compatible with the Mukai vector introduced in the next subsection). By \autoref{compdegree} (and the fact that the kernel bundle is stable, since $\Pic(S)=\bZ\cdot L$) the above set-theoretic description is equivalent to 
\[
W^2_d(S,L)=\bigsqcup_c W^2_d(S,L)_c,
\]
where
\[
W^2_d(S,L)_c:=\{(E,V^\vee) \;| \; E\in \cM(2,L,g+1-c),\; V^\vee\in \Gr(3,H^0(E)),\; \mathrm{deg}(\bigcap_{s\in V^\vee}Z_{cycle}(s))\geq c-d \}
\]
(note that the disjoint union is finite by \autoref{minimalc2}).

With this second description, $W^2_d(S,L)$ carries a natural scheme structure. Indeed, for a fixed $E\in\cM(2,L,g+1-c)$, one can show  that $\deg(\bigcap_{s\in V^\vee}Z_{cycle}(s))$ is an upper-semicontinuous function of $V^\vee\in\Gr(3,H^0(E))$; a relative construction makes $W^2_d(S,L)_c$ into a closed subset of the relative Grasmannian over $\cM(2,L,g+1-c)$ (see \cite[Section 2]{mor} for details).

We will consider these closed subsets with the reduced scheme structure\footnote{This is slightly different to \cite{mor}, where these families are pushed forward to $\Gr(3,H^0(L))$ via the injection induced by globalizing the maps $\mathrm{Gr}(3,H^0(E))\to \mathrm{Gr}(3,H^0(L))$, $V^\vee\mapsto \bigwedge ^2 V^\vee=V$.}.

The stratification of the \emph{Brill-Noether locus} $W^2_d(S,L)$ via $c_2(E)$ distinguishes rational maps by the length of the base locus (i.e. distinguishes the linear spaces we are projecting from by how secant to $S$ they are). For instance, for a $K3$ surface of genus $5$ we have 
\[
W^2_4(S,L)=W^2_4(S,L)_4\sqcup W^2_4(S,L)_5\cong W^2_4(S,L)_4 \sqcup S,
\] 
where $W^2_4(S,L)_4$ is the relative Grassmannian over $\cM(2,L,2)$ (see \cite[Theorem 1.2]{mor}).

More interestingly, if one fixes the second Chern class there may be different components according to the length of the base scheme of the sections of $V^\vee$. Namely
\[
W^2_d(S,L)_{c}=\bigcup_mW^2_d(S,L)_{c,m}
\]
where 
\[
W^2_d(S,L)_{c,m}=\{(E,V^\vee) \in W^2_d(S,L)_c\;|\; \mathrm{length}\left(\cap_{s\in V^\vee}Z(s)\right)\ge m\}.
\]

The length of $\bigcap_{s\in V^\vee}Z(s)$ can be seen as a measure of how tangent to $S$ is the linear space we are projecting from, whereas the difference $\deg\left(\bigcap_{s\in V^\vee}Z_{cycle}(s)\right)-\mathrm{length}\left(\bigcap_{s\in V^\vee}Z(s)\right)$ is a measure of how singular this tangency is. 

In the case where the finite subscheme $\bigcap_{s\in V^\vee}Z(s)$ is curvilinear (i.e. it lies on a smooth curve contained in $S$), we need the following lemma to study the relation between these two invariants:

 \begin{lem}\label{ideallemma}
     Let $E$ be a rank 2 vector bundle on $S$, and let $V^\vee\in\Gr(3,H^0(E))$. Assume that $\xi=\bigcap_{s\in V^\vee}Z(s)$ is curvilinear of length $m$ at a point $p$, and that the cycle $\bigcap_{s\in V^\vee}Z_{cycle}(s)$ has degree $f\geq m$ at $p$. Then there exist local coordinates $x,y$ around $p$ and $e_1,e_2$ local generators of $E$ at $p$, such that the image of the map
     \[
     V^\vee\otimes\cO_{S,p}\longrightarrow E_p
     \]
     is contained in $(x,y^m)\cdot e_1+(x^{\lceil\frac{f}{m}\rceil},x^{\lceil\frac{f-1}{m}\rceil}y,...,x^{\lceil\frac{1}{m}\rceil}y^{f-1},y^f)\cdot e_2$.
 \end{lem}
 \begin{proof}
     We choose the local coordinates $x,y$ such that $\cI_\xi=(x,y^m)$ locally at $p$. Given $e_1,e_2$ local generators of $E$ at $p$, let $\cI_j$ be the monomial ideal generated by monomials appearing in the $j$-th coordinate of sections of $V^\vee$; we may choose $e_1,e_2$ so that $x\in\cI_1$.

     In this framework, three general elements $s_1,s_2,s_3$ of $V^\vee$ (giving a basis) have local expressions
     \[
     s_i=\left((x-f_i(y))\cdot t_i(x,y),u_i(x,y)\right),
     \]
     where $f_i(y)\in (y^m)$ and $t_i(x,y)\in\cO_{S,p}^*$ (by the implicit function theorem).

     Of course, one has $u_i(f_i(y),y)\in (y^f)$. Actually, since $Z_{cycle}(s_1+\lambda s_2+\mu s_3)$ has degree $\geq f$ at $p$, one has
     \[
     u_1(f_{\lambda,\mu}(y),y)+\lambda u_2(f_{\lambda,\mu}(y),y)+\mu u_3(f_{\lambda,\mu}(y),y)\in (y^f)
     \]
     for every $\lambda,\mu$ and every expression $(x-f_1)t_1+\lambda(x-f_2)t_2+\mu(x-f_3)t_3=(x-f_{\lambda,\mu}(y))\cdot t_{\lambda,\mu}(x,y)$ with $f_{\lambda,\mu}\in (y^m)$ and $t_{\lambda,\mu}\in\cO_{S,p}^*$.

     This last condition implies that, after replacing $e_2$ by an appropriate combination of $e_1$ and $e_2$, one may assume that $x^ky^l\notin\cI_2$ as long as $f>mk+l$. Therefore, for the new choice of local generators one has $\cI_2\subset (x^{\lceil\frac{f}{m}\rceil},x^{\lceil\frac{f-1}{m}\rceil}y,...,x^{\lceil\frac{1}{m}\rceil}y^{f-1},y^f)$, which proves the assertion.
 \end{proof}

The first natural step for constructing Brill-Noether loci is to understand whether the cycle-theoretic and the scheme-theoretic intersection coincide. The following remark will be useful:

\begin{rem} \label{singu}
Consider the cohomology jump loci
\[
R_{c,m}:=\{(E,\xi)\in \cM(2,L,g+1-c)\times S^{[m]}\;|\; h^0(E\otimes \cI_\xi)\ge 3\}
\]
and the corresponding relative Grassmanian
\[
\mathcal G_{c,m}:=\{(E,\xi,V^\vee)\;|\; (E,\xi)\in R_,, \;V^\vee\in\Gr(3,H^0(E\otimes \cI_\xi))\}.
\]

Then the natural forgetful map $\psi_{c,m}:\mathcal G_{c,m}\longrightarrow W^2_{c-m}(S,L)_{c,m}$ is an isomorphism outside $\psi_{c,m}^{-1}(\mathrm{Im}(\psi_{c,m+1}))$. Indeed, there is a well defined inverse morphism 
\[
(E,V^\vee)\mapsto (E,V^\vee,\cap_{s\in V^\vee} Z(s))\in \mathcal G_{c,m}
\]
precisely outside the image of $\psi_{c,m+1}$. In particular:

\begin{itemize}
    \item If $\mathrm{dim}(\mathrm{Im}(\psi_{c,m+1}))<\mathrm{dim}(\mathrm{Im}(\psi_{c,m}))$ and $\mathcal{G}_{c,m}$ is smooth, then the singularities of $\psi_{c,m}(\mathcal{G}_{c,m})$ are contained in $\psi_{c,m+1}(\mathcal G_{c,m+1})$ and $\psi_{c,m}$ is birational.

    \item If $R_{c,m+1}$ is empty and $R_{c,m}$ is irreducible, then $\psi_{c,m}$ is an isomorphism onto its image and $\cG_{c,m}$ is isomorphic to an irreducible component of $W^2_{c-m}(S,L)$.
\end{itemize}
We will mainly deal with the minimal second Chern class (i.e. $c=\lfloor\frac{g+3}{2}\rfloor$); in that case, we will just write $R_m,\mathcal G_m,\psi_m$ for simplicity.
\end{rem}

To finish this subsection, let us point out an elementary lemma concerning the number of local generators of an ideal sheaf:

 \begin{lem}\label{numbergen}
     Let $X$ be a variety, $T\subset X$ a closed subscheme and $p\in T$ a closed point. Then $\dim_\bC\Ext^1_{\cO_X}(\cO_T,k(p))$ equals the minimal number of generators of $\cI_{T/X}$ at the point $p$.
 \end{lem}
 \begin{proof}
     By applying the functor $\Hom_{\cO_X}(-,k(p))$ to the short exact sequence
     \[
     0\longrightarrow \cI_{T/X}\longrightarrow \cO_X\longrightarrow \cO_T\longrightarrow 0
     \]
     we get an isomorphism $\Hom_{\cO_X}(\cI_{T/X},k(p))\cong \Ext^1_{\cO_X}(\cO_T,k(p))$. Now the claim follows from the adjunction $\Hom_{\cO_X}(\cI_{T/X},k(p))\cong\Hom_{k}(\cI_{T/X}\otimes k(p),k(p))$ and Nakayama's lemma.\qedhere
 \end{proof}

\subsection{Bridgeland stability on $K3$ surfaces}
\label{DerCatSec}
Let $(S,L)$ be a polarized $K3$ surface. In this subsection we review some basic facts on the \emph{$(\alpha,\beta)$-plane} of Bridgeland stability conditions associated to $L$. For a more detailed discussion, the reader may consult Bridgeland original work \cite{Bridgeland:Stab, Bridgeland:K3} or the recent survey \cite{MS2}. 

We will denote the Mukai vector of an object $E\in\Db(S)$ by
\[
v(E)=(v_0(E),v_1(E),v_2(E)):=\ch(E)\cdot\sqrt{\td(S)}=(\ch_0(E),\ch_1(E),\ch_0(E)+\ch_2(E))
\]
This gives an element of $\Lambda:=\bZ\oplus\Pic(S)\oplus\bZ$, which is a lattice equipped with the pairing
\[
\langle v,w\rangle:=-\chi(v,w)=v_1\cdot w_1-v_0\cdot w_2-v_2\cdot w_0
\]
Let us also recall that a class $\delta\in\Lambda$ is called \emph{spherical} if $\delta^2=-2$. 

Given a coherent sheaf $E\in\Coh(S)$, we define its \emph{slope} by $\mu_L(E)=\frac{L\cdot \ch_1(E)}{L^2\cdot \ch_0(E)}$ (with the convention $\mu_L(E)=+\infty$ if $\ch_0(E)=0$). This leads to the usual notion of $\mu_L$-stability.

For a real number $\beta\in\bR$, we consider the full subcategories of $\Coh(S)$
\begin{align*}
    \cT_{\beta}:=\{E\in\Coh(S)\mid\mu_L(Q)>\beta\text{ for all quotients $E\twoheadrightarrow Q$}\}\\
    \cF_{\beta}:=\{E\in\Coh(S)\mid\mu_L(E)\leq\beta\text{ for all subsheaves $F\hookrightarrow E$}\},
\end{align*}
forming a torsion pair. The corresponding tilt gives a bounded t-structure on $\Db(S)$ with heart
\[
\Coh^{\beta}(S):=\set{E\in\Db(S)\mid\cH^{-1}(E)\in\cF_{\beta}, \; \cH^{0}(E)\in\cT_{\beta},\; \cH^i(E)=0 \text{ for }i\neq0,-1}.
\]
Finally, for $(\alpha,\beta)\in\bR_{>0}\times\bR$, let  $Z_{\alpha,\beta}:K_0(\Db(S))\to\bC$ be the group homomorphism given by
\[
Z_{\alpha,\beta}(E):=-\left(v_2(E)-\beta L\cdot v_1(E)+(\frac{\beta^2}{2}-\frac{\alpha^2}{2})L^2\cdot v_0(E)\right)+i\left(L\cdot v_1(E)-\beta L^2\cdot v_0(E)\right).
\]
Now we have all the ingredients to state a fundamental result of Bridgeland:

\begin{thm}[\cite{Bridgeland:K3}]
    The pair $\sigma_{\alpha,\beta}:=(\Coh^\beta(S),Z_{\alpha,\beta})$ is a Bridgeland stability condition on $\Db(S)$ if, for every $\delta\in\Lambda$ spherical with $\Im Z_{\alpha,\beta}(\delta)=0$ and $\delta_0>0$, one has $\Re Z_{\alpha,\beta}(\delta)>0$.
\end{thm}

\begin{rem}
    The associated region of stability conditions is thus an upper-half plane, from which certain holes (corresponding to points $(\alpha,\beta)$ for which $Z_{\alpha,\beta}(\delta)=0$ for a spherical class $\delta$) have been removed. This can be visualized (for Picard rank 1) in \cite[Figure~1]{BB}.
\end{rem}

The fact that $\sigma_{\alpha,\beta}$ is a Bridgeland stability condition can be essentially summarized as follows:

\begin{enumerate}
    \item For every nonzero $E\in\Coh^\beta(S)$ one has the inequality $\Im Z_{\alpha,\beta}(E)\geq0$, and $\Re Z_{\alpha,\beta}(E)<0$ whenever $\Im Z_{\alpha,\beta}(E)=0$.
    
    \item Every object of $\Coh^\beta(S)$ admits a \emph{Harder-Narasimhan} (\emph{HN} for short) \emph{filtration}, with respect to the notion of stability defined by the \emph{tilt slope}
    \[
    \nu_{\alpha,\beta}(E):=\left\{
    \begin{array}{c l}
     \frac{-\Re Z_{\alpha,\beta}(E)}{\Im Z_{\alpha,\beta}(E)} & \Im Z_{\alpha,\beta}(E) >0\\
     +\infty & \Im Z_{\alpha,\beta}(E)=0\\
    \end{array}
    \right.
    \]
\end{enumerate}

It is worth mentioning that Bridgeland stability conditions can be deformed continuously (more precisely, their set can be given the structure of a complex manifold), thanks to a technical condition called the \emph{support property}. While we refer the reader to \cite{BayerDef} for details, here we content ourselves with one of the main applications: there is a locally finite wall and chamber structure, so that stability of objects remains unchanged along a chamber.

A \emph{numerical wall} for a class $v\in\Lambda$ is the region of $\bR_{>0}\times\bR$ determined by an equation of the form $\nu_{\alpha,\beta}(v)=\nu_{\alpha,\beta}(w)$, where $w\in\Lambda$ is a class non-proportional to $v$.
An \emph{actual wall} for $v$ is a subset of a numerical wall, at which the set of semistable objects of class $v$ changes.


The structure of the walls in the $(\alpha,\beta)$-plane admits the following description:

\begin{prop} \label{structurewalls}
Assume that $\Pic(S)=\bZ\cdot L$, and let $v\in\Lambda$ be any class.
\begin{enumerate}[{\rm (1)}]
 \item All numerical walls for $v$ are either semicircles centered on the $\beta$-axis or lines parallel to the $\alpha$-axis.
 \item If $v_0\neq0$, there is a unique vertical wall with equation $\beta=\frac{L\cdot v_1}{L^2\cdot v_0}$. If $v_0=0$, then there is no vertical wall.
 \item The curve $H_v:\nu_{\alpha,\beta}(v)=0$ intersects every semicircular wall at its top point.
 This curve is an hyperbola (if $v_0\neq0$ and $v^2>0$), a pair of lines (if $v_0\neq0$ and $v^2=0$), a parabola (if $v^2<0$) or a single vertical line (if $v_0=0$).
 \item If $v^2\geq0$, then all semicircular walls are strictly nested.
 \item If $v^2<0$, then all walls intersect at the unique point $(\alpha,\beta)$ such that $Z_{\alpha,\beta}(v)=0$.
\end{enumerate}
\end{prop}

One should note that numerical walls may cross holes of the $(\alpha,\beta)$-plane; such a hole may determine the end of an actual wall. Pictures illustrating \autoref{structurewalls} (including this phenomenon) may be found in \cite[Figures~1-3]{MS2}.

Using this structure of the walls, one can often understand how the stability of an object varies along certain regions of the $(\alpha,\beta)$-plane. This philosophy is specially well suited for Gieseker semistable sheaves, thanks to the following result:

\begin{prop}[\cite{Bridgeland:K3}]
Let $v\in\Lambda$ be a class with  $v_0>0$, and let $\beta<\frac{L\cdot v_1}{L^2\cdot v_0}$.
Then an object $F\in\Coh^\beta(S)$ of class $v(F)=v$ is $\sigma_{\alpha,\beta}$-semistable for every $\alpha\gg0$ if, and only if, $F$ is a Gieseker semistable sheaf (with respect to $L$).
\end{prop}

We close this subsection with an elementary trick that will be applied several times: 

\begin{lem}\label{minimalrank}
    Assume $\Pic(S)=\bZ\cdot L$. Let $\beta_0=\frac{a}{b}\in\bQ$ with $a,b$ coprime, and let $E\in\Coh^{\beta_0}(S)$ be an object with $(L\cdot v_1-\beta_0 L^2\cdot v_0)(E)=\frac{1}{b}L^2$. If $E$ is $\sigma_{\alpha_0,\beta_0}$-semistable for some $\alpha_0>0$, then $E$ is $\sigma_{\alpha,\beta_0}$-semistable for every $\alpha>0$.
\end{lem}

\subsection{Moduli spaces}
Moduli spaces of Gieseker stable sheaves on $K3$ surfaces are well understood thanks to the work of many authors after the pioneering work of Mukai (\cite{mukaisymp}). In the more general context of Bridgeland stability, the following result was established by Bayer-Macrì \cite{bayermacri}, building on previous work of Toda \cite{toda}:

\begin{thm}
Let $v\in\Lambda$ be a primitive class, and let $\sigma_{\alpha,\beta}$ be a Bridgeland stability condition lying in no actual wall for $v$. Then there exists a moduli space $\cM_{\sigma_{\alpha,\beta}}(v)$ of $\sigma_{\alpha,\beta}$-stable objects of class $v$, which is a smooth projective hyperkähler variety of dimension $v^2+2$ (in particular, it is nonempty if and only $v^2\geq -2$).
\end{thm}

We will be especially interested in the case of the Gieseker moduli space $\cM:=\cM(2,L,\frac{g-1}{2})$, for $(S,L)$ a polarized $K3$ surface of odd genus $g$ with $\Pic(S)=\bZ\cdot L$. This is a $K3$ surface, and results of Mukai (\cite{mvbnd}) show that:

\begin{enumerate}
    \item If $g\equiv 3\pmod{4}$, then $\cM$ is a fine moduli space and $\Pic(\cM)=\bZ\cdot\hat{L}$, where $\hat{L}$ is a polarization of genus $g$. Furthermore, the universal bundle $\cE$ on $S\times \cM$ induces a Fourier-Mukai equivalence $\Phi_\cE:\Db(S)\to\Db(\cM)$.

    \item If $g\equiv 1\pmod{4}$, then $\cM$ is a coarse moduli space and $\Pic(\cM)=\bZ\cdot\hat{L}$, where $\hat{L}$ is a polarization of genus $\frac{g+3}{4}$.
\end{enumerate}

In the second case, there exists an open analytic cover $\{U_i\}_{i\in I}$ of $\cM$ so that there exists a local universal sheaf $\cE_i$ on each $S\times U_i$. Furthermore, there is a (nontrivial) Brauer class $\alpha\in\mathrm{Br}(\cM)$ and isomorphisms $\varphi_{ij}:\cE_j|_{U_i\cap U_j}\to \cE_i|_{U_i\cap U_j}$ making $(\{\cE_i\},\{\varphi_{ij}\})$ into an $\alpha$-twisted sheaf $\cE$ on $S\times\cM$. This was proved by Caldararu (\cite{caldararu}). In analogy to the classical (untwisted) case of fine moduli spaces, Caldararu also showed that $\cE$ is the kernel of a Fourier-Mukai equivalence
\[
\Phi_\cE:\Db(S)\to \Db(\cM,\alpha),
\]
where $\Db(\cM,\alpha)$ is the bounded derived category of $\alpha$-twisted coherent sheaves on $\cM$.

Write $g=4n-3$. As proved independently by Sawon (\cite{sawon}) and Markushevich (\cite{markus}), one can use $\Phi_\cE$ to define an isomorphism
\[
S^{[n]}\overset{\cong}{\longrightarrow}\cM((0,\hat{L},k),\alpha),\;\;\;\;\;\; \xi\mapsto R^1\Phi_\cE(\cI_\xi)
\]
for some $k\in\bZ$, where $\cM((0,\hat{L},k),\alpha)$ denotes the moduli space of stable $\alpha$-twisted sheaves on $\cM$ with twisted Mukai vector $(0,\hat{L},k)$ (see \cite{twist} for the construction of the twisted Mukai vector depending on a B-field lift of $\alpha$, and \cite{yostwisted}
for the construction of moduli spaces of stable twisted sheaves).

By composing the isomorphism above with the natural map 
\[
\cM((0,\hat{L},k),\alpha)\longrightarrow \bP^n=\bP(H^0(\cM,\hat{L}))
\]
sending ($\alpha$-twisted) sheaves to their support, one gets a Lagrangian fibration $\pi$ of $S^{[n]}$.

Since Brauer classes on curves can be trivialized, any element in $\cM((0,\hat{L},k),\alpha)$ can be identified with an untwisted sheaf; this identification is not canonical, and hence $\cM((0,\hat{L},k),\alpha)$ is not globally isomorphic to the Gieseker moduli space $\cM(0,\hat{L},k)$ on $\cM$.

Nevertheless, by fixing a trivialization on any given curve $C\in |\hat{L}|=\bP(H^0(\cM,\hat{L}))$, it follows that the fiber $\pi^{-1}(C)$ is isomorphic to the compactified Jacobian $\overline{JC}$ parametrizing rank 1, torsion-free sheaves on $C$ of degree 0.

We can use these results to control the cohomology $h^1(S,E\otimes\cI_\xi)$ for $E\in\cM$ and $\xi\in S^{[n]}$:

\begin{lem} 
\label{coh913}
For the Lagrangian fibration $\pi:S^{[n]}\overset{\cong}{\longrightarrow}\cM((0,\hat{L},k),\alpha)\longrightarrow \bP^n=|\hat{L}|$, we have:
\begin{enumerate}[{\rm (1)}]
    \item $\pi$ sends $\xi\in S^{[n]}$ to the curve $C_\xi:=\{E\in\cM\, | \, h^1(E\otimes\cI_\xi)>0\}$
in the primitive linear system of $\cM$.

\item Given $\xi\in S^{[n]}$, let $L_\xi\in \overline{JC_\xi}$ denote the corresponding element under any identification $\pi^{-1}(C_\xi)\cong \overline{JC_\xi}$. Then for any $E\in C_\xi$, the number $h^1(S,E\otimes\cI_\xi)$ equals the minimal number of local generators of $L_\xi$ at the point $E$.

\noindent In particular, if $\mathrm{mult}_{E}(C_\xi)= m$, then $h^1(S, E\otimes\cI_\xi)\leq m$.
\end{enumerate}

\end{lem}
\begin{proof}
Since $h^2(E\otimes\cI_\xi)=0$ for every $E\in\cM$ and $\xi\in S^{[n]}$, we have an isomorphism 
\[
R^1\Phi_\cE(\cI_\xi)\otimes k(E)\overset{\cong}{\longrightarrow} H^1(E\otimes\cI_\xi)
\]
by cohomology and base change (which is still valid in the twisted setting, since it is a local property). Only the last assertion is not an immediate consequence of this isomorphism.

To complete the proof, note that all elements of $\overline{JC_\xi}$ are ideal sheaves, up to twist by a line bundle (see e.g. \cite[Proposition 3.2]{dsouza}). Since the curve $C_\xi$ is Gorenstein (it lies on a smooth surface), the minimal number of generators of an ideal sheaf at a point is bounded by the multiplicity of the curve at that point, which finishes the proof.
\end{proof}

\section{Genus 7}\label{sec:genus7}

Let $(S,L)$ be a polarized $K3$ surface of genus 7 (i.e. $L^2=12$). The goal of this section is to prove \autoref{thmA}.(1), hence we will assume $\Pic(S)=\bZ\cdot L$ throughout the rest of this section.
The upper bound $\irr_L(S)\leq 4$ was already obtained in \cite{mor}; let us recall the construction of rational dominant maps $S\dasharrow\bP^2$ of degree 4 given there.

We denote by $\cM=\cM(2,L,3)$ the fine moduli space of stable rank 2 bundles with $c_1(E)=L$ and $c_2(E)=5$. Any $E\in\cM$ has $h^0(E)=5$, hence $h^0(E\otimes\cI_p)\geq 3$ for any $p\in S$ (an equality holds, actually). Then the 3-dimensional subspace $V^\vee= H^0(E\otimes\cI_p)\subset H^0(E)$ induces by duality an element $V\in\Gr(3,H^0(L))$, such that $\deg(\varphi_V)\leq 4$ (by \autoref{compdegree}).

\subsection{Discarding rational maps of degree 3} A sufficient condition for having a map $S\dasharrow\bP^2$ of degree 3 is $h^1(E\otimes\cI_\xi)\geq 2$ for some $\xi\in S^{[2]}$ (i.e. $h^0(E\otimes\cI_\xi)\geq 3$). We will prove, by means of Bridgeland stability, that such a configuration is impossible.

Observe that $h^1(E\otimes\cI_\xi)=\hom(\cI_\xi, E^\vee[1])$. Furthermore, $\cI_\xi\in\Coh^\beta(S)$ for all $\beta<0$ (it is $\sigma_{\alpha,\frac{-1}{2}}$ and $\sigma_{\alpha,\frac{-1}{3}}$-stable for every $\alpha$ by \autoref{minimalrank}) and $E^\vee[1]\in\Coh^\beta(S)$ for all $\beta\geq\frac{-1}{2}$ (it is $\sigma_{\alpha,\beta}$-stable for all $\alpha$ as long as $\beta=\frac{-1}{3},\frac{-2}{5},\frac{-3}{7}$, by \autoref{minimalrank}).

The numerical wall $W=W(\cI_\xi,E^\vee[1])$ has equation $6(\beta^2+\alpha^2)+5\beta+1=0$. Since it intersects $\alpha=0$ at the points with $\beta$-coordinates $\frac{-1}{2}$ and $\frac{-1}{3}$, it is clear that both $\cI_\xi$ and $E^\vee[1]$ are semistable along $W$.

Now observe that any nonzero morphism $\cI_\xi\to E^\vee[1]$ defines a short exact sequence
\[
0\to F\to \cI_\xi\to E^\vee[1]\to 0
\]
in $\Coh^\beta(S)$ for all $\frac{-1}{2}\leq\beta<\frac{-1}{3}$ (indeed, such a morphism corresponds to a stable extension of sheaves $0\to E^\vee\to F\to  \cI_\xi\to 0$ with $\mu_L(F)=\frac{-1}{3}$). Moreover, $F$ is $\sigma_{\alpha,\beta}$-stable for all $\alpha$ as long as $\beta=\frac{-1}{2},\frac{-2}{5},\frac{-3}{8}$ (again by \autoref{minimalrank}). Since $\nu_{\alpha,\beta}(\cI_\xi)>\nu_{\alpha,\beta}(E^\vee[1])$ for $(\alpha,\beta)$ inside $W$, it follows that $\cI_\xi$ gets destabilized along $W$.

The situation can be visualized in the following picture. Note that $W$ crosses the hole in the $(\alpha,\beta)$-plane corresponding to the spherical bundle $G$ of Mukai vector $v(G)=(5,-2L,5)$. We will denote by $\sigma_0$ the Bridgeland stability condition $W\cap \{\beta=\frac{-3}{7}\}$.

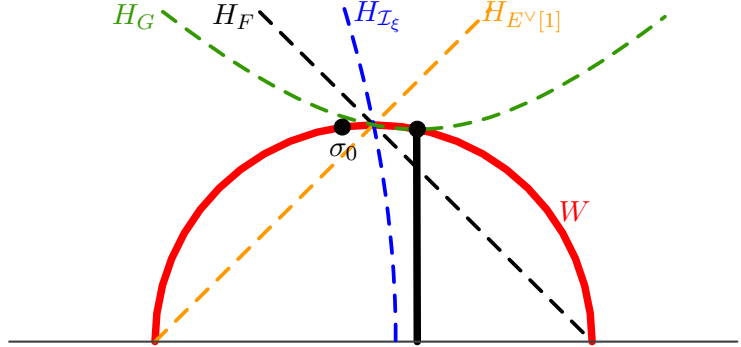
\begin{figure}[H]
\definecolor{uququq}{rgb}{0.25098039215686274,0.25098039215686274,0.25098039215686274}
\definecolor{uuuuuu}{rgb}{0.26666666666666666,0.26666666666666666,0.26666666666666666}
\definecolor{ffzzqq}{rgb}{1,0.6,0}
\definecolor{ttzzqq}{rgb}{0.2,0.6,0}
\definecolor{ffqqqq}{rgb}{1,0,0}
\definecolor{qqqqff}{rgb}{0,0,1}
\begin{tikzpicture}[line cap=round,line join=round,>=triangle 45,x=34.5cm,y=34.5cm]
\clip(-0.5899814447950904,-0.017083436252164914) rectangle (-0.24644106033089822,0.135);
\draw [line width=2.8pt] (-0.40022049484326594,0.08164965809277261)-- (-0.4,0);

\draw [line width=2.8pt,color=ffqqqq] plot[domain=0:3.141592653589793,variable=\t]({-0.4166666666666667+0.08333333333333331*cos(\t r)},{0.08333333333333331*sin(\t r)});
\draw [line width=1.4pt,dash pattern=on 6pt off 5pt,domain=-0.46:-0.333333333333] plot(\x,{-\x-0.33333333333333});
\draw [line width=1.4pt,dash pattern=on 6pt off 5pt,color=ffzzqq,domain=-0.5:-0.372] plot(\x,{(-1-2*\x)/-2});
\draw [line width=0.8pt,color=uququq,domain=-0.555555:-0.27777777] plot(\x,{(-0-0*\x)/1});
\draw [samples=50,domain=-0.153:0,line width=1.4pt,dash pattern=on 6pt off 5pt,color=qqqqff] plot ({5*0.08164965809277246*(-1-(\x)^2)/(1-(\x)^2)},{5*0.08164965809277246*(-2)*(\x)/(1-(\x)^2)});
\draw [samples=50,domain=-0.497:-0.305,line width=1.4pt,dash pattern=on 6pt off 5pt,color=ttzzqq] plot (\x,{sqrt((\x)^2+0.8*\x+0.166666666)});

\draw[color=ffqqqq] (-0.34,0.05) node {$W$};
\draw (-0.47,0.125) node {$H_F$};
\draw[color=ffzzqq] (-0.36,0.125) node {$H_{E^\vee[1]}$};
\draw[color=qqqqff] (-0.415,0.125) node {$H_{\mathcal{I}_\xi}$};
\draw[color=ttzzqq] (-0.508,0.125) node {$H_G$};
\draw (-0.428,0.073) node {$\sigma_0$};

\begin{scriptsize}
\draw[color=qqqqff] (0.8194876944657251,0.691774648870683) node {$ec5$};
\draw [fill] (-0.4,0.08164965809276782) circle (3pt);
\draw [fill] (-0.42857142857142855,0.08247860988423228) circle (3pt);
\end{scriptsize}
\end{tikzpicture}
\end{figure}

\begin{lem}\label{isomoduli} Let $F\in\cM(3,-L,2)$ and $E\in\cM$. Then:
\begin{enumerate}[{\rm (1)}]
 \item $F$ is destabilized along $W\cap\{\beta<\frac{-2}{5}\}$, and $E^\vee[1]$ is destabilized along $W\cap\{\beta>\frac{-2}{5}\}$.

 \item There is an isomorphism $\varphi:\cM\overset{\cong}{\longrightarrow}\cM(3,-L,2)$.
\end{enumerate} 
\end{lem}
\begin{proof}
    Let $E\in\cM$. Note that $G\in\Coh^\beta(S)$ for all $\beta<\frac{-2}{5}$, and it is $\sigma_{\alpha,\beta}$-stable for all $\alpha$ if $\beta=\frac{-1}{2},\frac{-3}{7}$; in particular, $G$ and $E^\vee[1]$ are $\sigma_0$-stable of the same slope. It follows that $\hom(E^\vee[1],G)=0=\hom(G,E^\vee[1])$, which gives $\ext^1(E^\vee[1],G)=\langle v(E^\vee[1],v(G)\rangle=1$.

    Therefore there is a (unique) nontrivial extension
    \begin{equation}\label{extwall7}
        0\to G\to \varphi(E)\to E^\vee[1]\to 0,
    \end{equation}
    and $\varphi(E)$ must be $\sigma$-semistable for $\sigma\in\{\beta=\frac{-3}{7}\}$ above $\sigma_0$. Otherwise, for a $\sigma$-destabilizing subobject $R\subset \varphi(E)$ the composition $R\hookrightarrow \varphi(E)\twoheadrightarrow E^\vee[1]$ would be nonzero (since $G$ is $\sigma$-stable), and $\sigma_0$-stability of $E^\vee[1]$ would force $R\cong E^\vee[1]$, which splits the extension \eqref{extwall7}.

    By $\sigma$-stability of $\varphi(E)$ it follows that $\varphi(E)\in\cM(3,-L,2)$, and since this construction can be carried out in families\footnote{Let $\cE$ be a universal bundle on $S\times\cM$, and let $p,q$ denote the two projections of $S\times\cM$ onto $S$ and $\cM$. By cohomology and base change, $q_*(p^*G\otimes\cE)$ is a line bundle whose fiber at $E\in\cM$ is canonically identified with $\Ext^1(E^\vee[1],G)$. For an appropriate choice of $\cE$, we may assume $q_*(p^*G\otimes\cE)\cong\cO_\cM$; then we obtain a unique (up to constant) morphism $\cE^\vee \to p^*G$, which globalizes the construction above.}, the assignment $\varphi$ is a morphism of schemes. Furthermore, the sequence \eqref{extwall7} destabilizes $\varphi(E)$ along $W\cap\{\beta<\frac{-2}{5}\}$: $G$ and $E^\vee[1]$ are the HN factors of $\varphi(E)$ for stability conditions along $\beta=\frac{-3}{7}$ that are below $\sigma_0$.

    By rotating the distinguished triangle, we obtain a short exact sequence
    \[
    0\to \varphi(E)\to E^\vee[1]\to G[1]\to 0
    \]
    in $\Coh^\beta(S)$ ($\frac{-2}{5}\leq\beta<\frac{-1}{3}$), which destabilizes $E^\vee[1]$ along $W\cap \{\beta>\frac{-2}{5}\}$. The HN factors of $E^\vee[1]$ for stability conditions just below $W\cap\{\beta=\frac{-3}{8}\}$ will be precisely $\varphi(E)$ and $G[1]$.

    On the other hand, starting with $F\in\cM(3,-L,2)$, one can construct an element $\psi(F)\in\cM$ arising as the (unique) nontrivial extension
    \[
    0\to F\to \psi(F)^\vee[1]\to G[1]\to 0
    \]
    (for the uniqueness of such an extension, we use again the stability of $F$ and $G[1]$ along $\beta=\frac{-3}{8}$). By rotating the triangle we can destabilize $F$ along $W\cap\{\beta<\frac{-2}{5}\}$.

    The morphisms $\varphi$ and $\psi$ are inverse to each other, as easily follows from the uniqueness of the HN filtrations of $E^\vee[1]$, $\varphi(E)$, $F$ and $\psi(F)^\vee[1]$ just below $W$. This concludes the proof.
\end{proof}

Now we can prove the desired cohomological bound:

\begin{prop}\label{cohboundgenus7}
    For all $E\in\cM$ and $\xi\in S^{[2]}$, the inequality $h^1(E\otimes\cI_\xi)\leq 1$ holds.
\end{prop}
\begin{proof}
    Assume that $h^1(E\otimes\cI_\xi)\geq 1$, i.e., there exists a short exact sequence
    \[
    0\to F\to \cI_\xi\to E^\vee[1]\to 0
    \]
    destabilizing $\cI_\xi$ along $W$ (with $F\in\cM(3,-L,2)$). Consider the short exact sequence 
    \begin{equation}\label{destF}
        0\to G\to F\to \psi(F)^\vee[1]\to 0
    \end{equation}
    destabilizing $F$ along $W\cap\{\beta<\frac{-2}{5}\}$, constructed in the proof of \autoref{isomoduli}. 
    Then the HN factors of $\cI_\xi$ (for stability conditions just below $\sigma_0$) are $G$ and $T$, where $T$ is the extension of $E^\vee[1]$ by $\psi(F)^\vee[1]$ obtained as the image of $\cI_\xi$ under the induced map
    \[
    \Ext^1(E^\vee[1],F)\longrightarrow\Ext^1(E^\vee[1],\psi(F)^\vee[1])
    \]

    If $E\neq \psi(F)$, of course we have $T=E^\vee[1]\oplus\psi(F)^\vee[1]$; otherwise, $T$ is a nontrivial extension. Indeed, if $E= \psi(F)$ then the above map of Ext-groups is an isomorphism (as follows from the long exact sequence obtained by applying the functor $\Hom(E^\vee[1],-)$ to \eqref{destF}).

    In any case, since $\Hom(G,E^\vee[1])=0$ (recall that $G$ are $E^\vee[1]$ are $\sigma_0$-stable of the same slope), we obtain that $\hom(\cI_\xi, E^\vee[1])=\hom(T,E^\vee[1])=1$.
\end{proof}

Recall that the kernel bundle of any degree 3 map has $c_2=5$ (\autoref{minimalc2}). Hence the following lemma, together with \autoref{cohboundgenus7}, discards the existence of degree 3 maps and concludes the proof of the equality $\mathrm{irr}_L(S)=4$:
	
 \begin{lem}
	Let $E \in \mathcal M$. Then there exists no $V^\vee \in \mathrm{Gr}(3,H^0(E))$ such that $\bigcap_{s\in  V^\vee}Z_{cycle}(s)=2p$ and $\bigcap_{s \in V^\vee}Z(s)=\{p\}$ (i.e. the intersection as zero cycles has degree $2$, whereas the schematic intersection of the zero loci of the sections is a reduced point).
	\label{vvv}
	\end{lem}
\begin{proof}
We argue by contradiction. Suppose that there exists $V^\vee \in \mathrm{Gr}(3,H^0(E))$ such that $\bigcap_{s\in  V^\vee}Z_{cycle}(s)=2p$ and $\bigcap_{s \in V^\vee}Z(s)=\{p\}$. According to \autoref{ideallemma}, we can choose local coordinates and local generators $e_1,e_2$ of $E$ at $p$, such that the subsheaf $E_1 \subset E$ defined by ${E_1}|_{S\setminus p}=E|_{S\setminus p}$ and $E_1=m_p e_1+m_p^2e_2$ in a neighborhood of $p$ has $V^\vee \subset H^0(E_1)$, in particular we get $h^0(E_1)\ge 3$. Since the quotient $E/E_1$ is a $4$-dimensional vector space, we get $\chi(E_1)=1$ and hence $h^1(E_1)=2$. Let us consider the double extension
	\[
	\begin{tikzcd}
		0 \arrow{r} & \mathcal O_S^{\oplus 2} \arrow{r} & F \arrow{r} & E_1 \arrow{r} & 0.
	\end{tikzcd}\]
It is easy to compute that $h^0(E_1 \otimes k(p))=5$ ($E_1$ has $xe_1$, $ye_1$, $x^2e_2$, $xye_2$, $y^2e_2$ as local generators), hence $F$ cannot be locally free at $p$. We get $c_2(F^{**})<c_2(F)=c_2(E_1)=9$. It suffices to prove that $F^{**}$ must be stable, since there exists no stable rank $4$ vector bundle with these invariants. To this end, consider the commutative diagram with exact rows
\[
\begin{tikzcd} 0 \arrow{r} & E_1^*=E^\vee \arrow{r} \arrow{d} & V\otimes \mathcal O_S \arrow{r} \arrow{d} & L \arrow{d}{=} \\
    0 \arrow{r} & F^* \arrow{r}\arrow{d} & W \otimes \mathcal O_S \arrow{r} \arrow{d}& L \arrow{d} \\
   0\arrow{r} & \mathcal O_S^{\oplus 2} \arrow{r}{=} & \mathcal O_S^{\oplus 2}\arrow{r}& 0\end{tikzcd}
\]
where $V^\vee\in \mathrm{Gr}(3,H^0(E_1))$. The map $W\to H^0(L)$ is injective; indeed, $V \to H^0(L)$ is injective (by stability of $E$) and the extension defining $F$ is nonsplit (none of the maps $\mathcal O_S\to F$ in the definition of $F$ split). It follows that $H^0(F^*)=0$, which is enough to conclude the stability of the kernel bundle $F^*$ since $\Pic(S)=\bZ\cdot L$ (if $T$ is a vector bundle destabilizing $F^*$, then $T$ must be of the form $\cO_S^{\oplus\rk T}$ since $T\hookrightarrow W\otimes\cO_S$, contradiction).
\end{proof}

\subsection{Determining $W^2_4(S,L)$}

In the notations of \autoref{singu}, $R_2$ is empty. Hence $\psi_1:S \times \cM \cong \cG_1\to W^2_4(S,L)_{5,1}$ is an isomorphism.
The following lemma shows that this is the unique component of the Brill-Noether locus $W^2_4(S,L)$:
\begin{lem}
    If $S \dashrightarrow \mathbb P^2$ is a rational map of degree $4$, then the kernel bundle $E$ has $c_2(E)=5$.
\end{lem}
\begin{proof}

By \autoref{minimalc2} the only other possibility is $c_2(E)=6$ and $V^\vee \in \mathrm{Gr}(3,H^0(E))$ with $\bigcap_{s\in  V^\vee}Z_{cycle}(s)$ of degree $2$.

If $\bigcap_{s\in  V^\vee}Z(s)$ is a reduced point, then by \autoref{ideallemma} in an appropriate trivialization of $E$ near $p$ all sections $s\in V^\vee$ are of the form $s=ae_1+be_2$, with $a\in m_p\setminus m_p^2$ and $b\in m_p^2$. Since $V=\bigwedge ^2 V^\vee$, it follows that the base ideal $\cI$ satisfies $\cI\subset \cI_p^3$; therefore, $\mathrm{colength}(\cI)\geq 7$ (because $\cI$ is generated by three elements). This implies that $c_2(E)=L^2-\mathrm{colength}(\cI)\leq5$. 

Hence, if $c_2(E)=6$, we must have  $h^0(E \otimes \mathcal I_\xi)\ge 3$  for some $\xi \in S^{[2]}$. Consider now the subsheaf $E_1  \subset E$ generated by $E \otimes \mathcal I_\xi$ and $\mathrm{Im}(\mathcal O_S \otimes H^0(E) \to E)$. The quotient $E_1/E\otimes \mathcal I_\xi$ is of length at most $2$, since it is a principal $\mathcal{O}_\xi$-submodule in $E/E \otimes \mathcal I_\xi$ (it is generated by a section in $H^0(E)\setminus H^0(E \otimes \mathcal I_\xi)$). This implies $\chi(E_1)\leq \chi(E \otimes \mathcal I_\xi)+2=2$, and since $h^0(E_1)=4$, we obtain $h^1(E_1)\ge 2$. Now the nontrivial extension
\[
\begin{tikzcd}
0 \arrow{r} & \mathcal{O}_S^{\oplus h^1(E_1)} \arrow{r} & F \arrow{r} & E_1 \arrow{r} & 0.
\end{tikzcd}
\]
gives a stable torsion free sheaf $F$ with Mukai vector $v(F)=(h^1(E_1)+2,L,2)$, which yields a contradiction if $h^1(E_1)\geq2$.
\end{proof}

\subsection{An isomorphism of Hilbert schemes}\label{Hilbschemes} Before we move to higher genera, let us use the framework provided by Bridgeland stability to establish an isomorphism $S^{[2]}\cong\cM^{[2]}$ (already predicted by Yoshioka \cite{yos}). The arguments here will be useful in genus 11.

\begin{prop}\label{isohilb2}
The Hilbert schemes $S^{[2]}$ and $\cM^{[2]}$ are isomorphic.
\end{prop}
\begin{proof}
Let us define a morphism $\gamma:\cM^{[2]}\longrightarrow S^{[2]}$ as follows. 

Given a reduced element of $\cM^{[2]}$, supported at two distinct $E_1,E_2\in\cM$, we consider the object $M$ fitting into an extension
\[
0\to G\to M\to E_1^\vee[1]\oplus E_2^\vee[1]\to 0
\]
such that the induced elements of $\Ext^1(E_i^\vee[1],G)$ are nonzero. Due to the action of automorphisms of $E_1^\vee[1]\oplus E_2^\vee[1]$, such an object $M$ is unique up to isomorphism.

On the other hand, given a nonreduced element of $\cM^{[2]}$ supported at $E\in\cM$ with tangent direction $T\in\bP(\Ext^1(E,E))$, we consider the corresponding $\widetilde{T}:=T^\vee[1]\in\bP(\Ext^1(E^\vee[1],E^\vee[1]))$ and take the object $M$ fitting into 
\[
0\to G\to M\to \widetilde{T}\to 0,
\]
such that the induced extension in $\Ext^1(E^\vee[1],G)$ is nonzero. Again, such an object $M$ is unique up to isomorphism (due to the action of automorphisms of $\widetilde{T}$).

Let us assume for a moment that, in both cases, the object $M$ is of the form $M=\cI_\xi$ for some $\xi\in S^{[2]}$ (we will prove it below). Then this rule defines the morphism $\gamma$, which is injective on closed points (indeed, given $M=\cI_\xi$ in the image of $\gamma$, the corresponding element of $\cM^{[2]}$ can be recovered from the last HN factor of $M$ just below $\sigma_0$). Therefore $S^{[2]}$ equals the image of $\gamma$ (both $\cM^{[2]}$ and $S^{[2]}$ are projective of the same dimension). Since $\gamma$ is bijective on closed points and $S^{[2]}$ is normal, it follows that $\gamma$ is an isomorphism.

Hence it only remains to prove that the objects $M$ constructed above are of the form $M=\cI_\xi$, $\xi\in S^{[2]}$. Since there are no walls above $W$ for the Mukai vector $(1,0,-1)$, it suffices to check that $M$ is $\sigma'$-semistable for some stability condition $\sigma'$ above $W$. We will assume that $M$ comes from a nonreduced element of $\cM^{[2]}$ (the reduced case being similar but simpler).

Let $\sigma'$ be close enough to $\sigma_0$ (above $W$), and assume that $M$ is not $\sigma'$-semistable. Since $M$ is $\sigma_0$-semistable, any destabilizing subobject $A\subset M$ is $\sigma_0$-semistable of the same slope as $M$.

The composition $G\hookrightarrow M\twoheadrightarrow M/A$ is either zero or injective, since $G$ is $\sigma_0$-stable of the same slope as $M/A$. In the first case, we have $G\hookrightarrow A$ and hence $A/G\hookrightarrow \widetilde{T}$. Consider the composition $A/G\hookrightarrow \widetilde{T}\twoheadrightarrow E^\vee[1]$. By $\sigma_0$-stability of $E^\vee[1]$, this morphism is either zero or surjective.  If $A/G\to E^\vee[1]$ is zero then either $A/G=0$ or $A/G\cong E^\vee[1]=\ker(\widetilde{T}\twoheadrightarrow E^\vee[1])$ (hence $\nu_{\sigma'}(A)<\nu_{\sigma'}(M)$, contradiction); whereas if $A/G\hookrightarrow \widetilde{T}\twoheadrightarrow E^\vee[1]$ is surjective, then either $A/G=\widetilde{T}$ (i.e. $A=M$, contradiction) or the extension $0\to E^\vee[1]\to \widetilde{T}\to E^\vee[1]\to 0$ splits.

On the other hand, if $G\hookrightarrow M/A$ is injective, one easily checks that $A\hookrightarrow \widetilde{T}$. Now we consider the composition $A\hookrightarrow \widetilde{T}\twoheadrightarrow E^\vee[1]$, which again will be either zero or surjective. If it is surjective, then either $A\cong \widetilde{T}$ (which splits $0\to G\to M\to \widetilde{T}\to 0$) or $A\cong E^\vee[1]$ (which splits $0\to E^\vee[1]\to \widetilde{T}\to E^\vee[1]\to 0$). If $A\hookrightarrow \widetilde{T}\twoheadrightarrow E^\vee[1]$ is zero, then $A\cong E^\vee[1]=\ker(\widetilde{T}\twoheadrightarrow E^\vee[1])$ and the induced element in $\Ext^1(E^\vee[1],G)$ is zero, which gives again a contradiction.
\end{proof}

\begin{rem}
    Let $\cE$ denote the universal vector bundle on $S\times\cM$, and $\Phi_\cE:\Db(\cM)\to \Db(S)$ the corresponding Fourier-Mukai equivalence. The isomorphism $\gamma$ in \autoref{isohilb2} sends $T\in\cM^{[2]}$ to the subscheme $\gamma(T)\in S^{[2]}$ whose ideal sheaf sits in a distinguished triangle $0\to G\to \cI_{\gamma(T)}\to \Phi_\cE(\cO_T)^\vee[1]\to 0$.
\end{rem}

\section{Genus 9}\label{sec:genus9}

Let $(S,L)$ be a polarized $K3$ surface of genus 9 (i.e. $L^2=16$) with $\Pic(S)=\bZ\cdot L$. Throughout this section, $\cM$ will denote the Gieseker moduli space $\cM(2,L,4)$.

We have $\irr_L(S)\leq 5$, as proved in \cite{mor} by considering vector bundles $E\in\cM$. Indeed, $h^0(E\otimes\cI_p)=4$ for every $p\in S$ and, since $c_2(E)=6$, any 3-dimensional vector space $V^\vee\subset H^0(E\otimes\cI_p)$ produces a map $\varphi_V:S\dasharrow\bP^2$ of degree $\leq 5$.

\subsection{Constructing maps of degree 4 (I)}\label{deg4genus9}
We start by constructing maps of degree 4, by considering pairs $(E,\xi)\in\cM\times S^{[2]}$ such that $h^1(E\otimes\cI_\xi)\geq 1$ (i.e. $h^0(E\otimes\cI_\xi)\geq 3$). 

To this end, we consider the numerical wall $W$ in the $(\alpha,\beta)$-plane defined by the Mukai vectors $(-2,L,-4)$ and $(1,0,-1)$. This wall intersects $\alpha=0$ at the points with $\beta$-coordinates $\frac{-1}{2}$ and $\frac{-1}{4}$; in particular, for every $\xi\in S^{[2]}$ the ideal sheaf $\cI_\xi$ is semistable along $W$, as an immediate application of \autoref{minimalrank}. Furthermore, $W$ crosses the hole in the $(\alpha,\beta)$-plane corresponding to the spherical bundle $F$ with Mukai vector $v(F)=(3,-L,3)$. Then:

\begin{enumerate}
    \item If $(E,\xi)\in\cM\times S^{[2]}$ satisfies $\hom(\cI_\xi, E^\vee[1])=h^1(E\otimes\cI_\xi)\geq 1$, then any nonzero map $\cI_\xi \to E^\vee[1]$ defines a destabilizing short exact sequence
    \[
    0\to F\to \cI_\xi \to E^\vee[1]\to 0
    \]
    for $\cI_\xi$ along $W\cap \{\beta<\frac{-1}{3}\}$. Moreover, since both $F$ and $E^\vee[1]$ are $\sigma_0$-stable for $\sigma_0=W\cap \{\beta=\frac{-2}{5}\}$ (this follows again from \autoref{minimalrank}), the short exact sequence above is the HN filtration of $\cI_\xi$ just below $\sigma_0$. In particular, the morphism $\cI_\xi\to E^\vee[1]$ is unique (up to constant).

    \item Conversely, it is easy to check that any nontrivial extension $0\to F\to G\to E^\vee[1]\to 0$ of $E^\vee[1]$ by $F$ must be semistable above $W$. Therefore $G=\cI_\xi$ for some $\xi\in S^{[2]}$, since there are no actual walls for the Mukai vector $(1,0,-1)$ above $W$.

    \noindent According to the vanishings $\hom(F,E^\vee[1])=0=\hom(E^\vee[1],F)$ (recall that $F$ and $E^\vee[1]$ are $\sigma_0$-stable of the same slope), we have
    \[
    \ext^1(E^\vee[1],F)=\langle v(E^\vee[1]),v(F)\rangle = 2.
    \]
    In other words, for a fixed $E$ the isomorphism classes of such nontrivial extensions are parametrized by a $\bP^1$.
\end{enumerate}

We can summarize this discussion as follows:

\begin{prop}\label{999}
    The locus $R_2=\{(E,\xi)\in \cM\times S^{[2]}\;|\;h^1(E\otimes\cI_\xi)\geq 1\}$ is a $\bP^1$-bundle over $\cM$, with a closed immersion $R_2\hookrightarrow S^{[2]}$ defined by the second projection.
\end{prop}

Let us also remark the following consequence of the above discussion:

\begin{lem}
    Let $F$ be the spherical vector bundle with Mukai vector $v(F)=(3,-L,3)$. Then for any $\xi \in S^{[2]}$, we have $\hom(F, \cI_\xi)\le 1.$
    \label{est}
\end{lem}
\begin{proof}
    This follows from the uniqueness of the HN filtration of $\cI_\xi$ just below $\sigma_0$.
\end{proof}

 \subsection{Discarding degree $3$ maps}
According to \autoref{compdegree} and \autoref{minimalc2}, a map $S\dasharrow\bP^2$ of degree 3 is equivalent to a pair $(E,V^\vee)$ with $E\in \cM$ (i.e. $c_2(E)=6$, $h^0(E)=6$) and $\deg\left(\bigcap_{s\in V^\vee}Z_{cycle}(s)\right)=3$. The base ideal $\cI$ of such a map has
\[
\colength(\cI)=L^2-c_2(E)=10.
\]
Recall that $\Pic(\cM)=\bZ\cdot\hat L$, where $\hat L$ is a degree $4$ polarization.

\begin{lem}
There is an isomorphism
\begin{align*}
R_3=\{ (E,\xi)\in \cM\times S^{[3]}\;|\; h^1(E\otimes \cI_\xi)\ge 3\}\cong \\
\cong \{ (C,E,A) \;|\; C\in &|\hat L|,E\in C, A\in \overline {JC}\text{ with }h^0(C,A\otimes k(E))\geq3\}. 
\end{align*}
In particular, if $(S,L)$ is very general among $K3$ surfaces of genus $9$, we have $h^1(E\otimes \cI_\xi)\le 2$ for every $(E,\xi)\in \cM\times S^{[3]}$.
\label{est1}
\end{lem}
\begin{proof}
The isomorphism is just a reformulation of \autoref{coh913}. For the second part, observe that the non-emptiness of $R_3$ implies that some plane quartic of the linear system $|\hat{L}|$ has a triple point, hence it is rational. If $(S,L)$ is very general, then $(\cM,\hat L)$ must be very general among quartic surfaces; but in that case a theorem of Chen \cite{chen} ensures that all the rational curves in the linear system $|\hat{L}|$ must be nodal, which gives a contradiction.
\end{proof}

\begin{rem}
   On the other hand, if $\Pic(S)=\bZ\cdot L$ and there is a curve with a triple point in $|\hat L|$, then the above discussion immediately yields $W^2_3(S,L)\neq \emptyset$. It would be interesting to know whether such $K3$ surfaces may exist.
\end{rem}

We are left to rule out the situation where the schematic intersection $\xi:=\bigcap_{s\in V^\vee}Z(s) $ is of length $\le 2$ and $\bigcap_{s\in V^\vee}Z_{cycle}(s)$ is of degree $3$. First let us slightly simplify the situation:

\begin{lem}
Under these assumptions, $\xi$ is of length $2$ and reduced.
\label{dajee}
\end{lem}
\begin{proof}
If $\xi$ is a reduced point $p$, then by \autoref{ideallemma} we have $s_p\in m_p\cdot e_1+m_p^3\cdot e_2$ for the local equations of $s\in V^\vee$ in an appropriate trivialization $E_p\cong e_1\cdot \cO_{S,p}\oplus e_2\cdot \cO_{S,p}$.  In view of the isomorphism $V\cong\bigwedge^2V^\vee$, it follows that the base ideal $\cI$ satisfies $\cI\otimes \cO_{S,p}\subset m_p^4$ locally at $p$; this implies that the degree of the map is $\leq L^2-16=0$, contradiction.

Now assume that $\xi$ is not reduced, locally of the form $(x,y^2)$ at $p=\mathrm{Supp}(\xi)$. In that case, by \autoref{ideallemma} we have $s_p\in(x,y^2)\cdot e_1+(x^2,xy,y^{3})\cdot e_2$ for all $s\in V^\vee$, in an appropriate trivialization $E_p\cong e_1\cdot \cO_{S,p}\oplus e_2\cdot \cO_{S,p}$ of $E$ near $p$. By the isomorphism $V\cong\bigwedge^2V^\vee$, it follows that the base ideal $\cI$ satisfies $\cI\otimes \hat\cO_{S,p}\subset (x^3,x^2y, xy^{3},y^{5})$ locally at $p$. Furthermore, since $\cI$ has at most three local generators at $p$, we obtain $\mathrm{colength}(\cI)\ge 11$, contradiction.
\end{proof}

 We are now ready to prove:
 
 \begin{cor}
 If $(S,L)$ is a very general polarized $K3$ surface of genus $9$, then $\mathrm{irr}_L(S)=4$.
 \label{ccc}
 \end{cor}
\begin{proof}
In view of \autoref{est1} and \autoref{dajee}, the only possibility for a degree 3 map is that $\bigcap_{s\in V^\vee} Z(s)$ is reduced, consisting of two distinct points $p,q$, and $\bigcap_{s\in V^\vee} Z_{cycle}(s)=2p+q$. A new application of \autoref{ideallemma} yields $V\cong \bigwedge^2 V^\vee\subset H^0(L\otimes \cI_p^3\otimes \cI_q^2)$. We deduce $h^1(L\otimes \cI_p^3\otimes \cI_q^2)=2$, hence we may consider the stable extension \[
\begin{tikzcd}
    0 \arrow{r} & \cO_S^{\oplus 2} \arrow{r} & G \arrow{r} & L \otimes \cI_p^3\otimes \cI_q^2 \arrow{r} & 0.
\end{tikzcd}
\]
It is easy to see that $G$ is not locally free at $p$ ($G$ has rank 3, whereas $L \otimes \cI_p^3\otimes \cI_q^2$ has four local generators). Hence $G\subset G^{**}$ is a nontrivial subsheaf, in particular $c_2(G^{**})<c_2(G)=9$, whereas on the other hand we have $c_2(G^{**})\ge 8$ since $G^{**}$ is stable  (stability follows from the fact that $G^*$ is a kernel bundle, with a proof similar to the stability of $F^*$ in the proof of \autoref{vvv}). This forces $G^{**}$ to be the spherical vector bundle with $v(G^{**})=(3,L,3)$. 

Moreover, $G^{**}$ fits in a commutative diagram with exact rows and columns 
\[
\begin{tikzcd}
&  & 0 \arrow{d} & 0 \arrow{d} & \\
0 \arrow{r} & \cO_S^{\oplus 2}\arrow{d}{=} \arrow{r} & G \arrow{r}\arrow{d} & L  \otimes \cI_p^3\otimes \cI_q^2 \arrow{r}\arrow{d} & 0\\
    0 \arrow{r} & \cO_S^{\oplus 2} \arrow{r} & G^{**} \arrow{r} \arrow{d} & L  \otimes \mathcal J\arrow{r} \arrow{d} & 0 \\
    & & k(p)\arrow{d}\arrow{r}{=} & k(p)\arrow{d}\\
    &   & 0 & 0 & 
\end{tikzcd}
\]
where $\mathcal J$ is an ideal of colength $8$ such that $\mathcal J/\mathcal I_p^3\otimes \mathcal I_q^2\cong k(p)$, hence it must have $3$ local generators at both $p$ and $q$. This implies that the map $\mathcal O_S^{\oplus 2}\to G^{**}$ vanishes at $p$ and $q$, hence $h^0(G^{**}\otimes \mathcal I_{p,q})\ge 2$ which is a contradiction to \autoref{est}. The proof is complete.
\end{proof}

\subsection{Constructing maps of degree 4 (II)}
 Here we are going to present another construction of degree $4$ maps, eventually leading to an irreducible component of $W^2_4(S,L)$ which is different from the $\bP^1$-bundle $R_2$ over $\cM$ that we described in \autoref{deg4genus9}.
 
 The idea is the following: given a nonreduced $\xi\in S^{[2]}$ (supported at $p\in S$) and $E\in \cM$ such that $h^1(E\otimes \cI_\xi)=1$, we want to find suitable local generators $e_1,e_2$ of $E$ around $p$ such that the subsheaf $\tilde{E}\subset E$ locally defined as $m_pe_1+m_p^2e_2$ (and coinciding with $E$ elsewhere) has $3$ global sections. One would expect the base locus of $\bigwedge ^2 H^0(E\otimes \cI_\xi)$ to be $\cI_\xi^2$ at $p$, whereas the base locus of $\bigwedge^2 H^0(\tilde{E})$ is locally  contained in $m_p^3$. Hence the two constructions should produce different rational maps, as we will indeed check.

Let us denote by $\pi_\cM:\cM\times S^{[2]}\to \cM$ the canonical projection, by $S^{[2]}_{nred}\subset S^{[2]}$ the divisor of nonreduced subschemes and by $\pi_S:\cM\times S^{[2]}_{nred}\to S$ the map sending a nonreduced subscheme to its support. Finally, for any $E \in \cM$ let us denote by $\mathbb P^1_E$ the immersion $\mathbb P^1_E=\mathbb P\left(\mathrm{Ext}^1(E^\vee[1],F)\right) \hookrightarrow S^{[2]}$ (recall \autoref{999}).

We need a preliminary lemma on the divisor
\[
R_2=\{(E,\xi) \in \cM\times S^{[2]} \;|\;h^1(E\otimes \cI_\xi)\ge 1\}=\bigcup_{E\in \cM} \mathbb P^1_E\hookrightarrow S^{[2]}
\]
constructed in \autoref{999}. 

\begin{lem}
The intersection
\[
\{(E,\xi) \in \cM\times S^{[2]} \; | \; \xi \textrm{ is nonreduced and } h^1(E\otimes \cI_\xi)=1\}
\]
of $R_2$ with the divisor $S^{[2]}_{nred}$ of nonreduced subschemes is $2$-dimensional (and not empty).
Furthermore, there is an irreducible component $W_S\subset R_2\cap S^{[2]}_{nred}$ such that the two projections $\pi_\cM:W_S\to \cM$ and $\pi_S:W_S\to S$ are surjective.
\label{reg}
\end{lem}
\begin{proof}
First let us show that $R_2 \not\subset S_{nred}^{[2]}$. Assume the contrary. If $\pi_S(\bP^1_E)$ is a point for every $E\in\cM$, then we get an injective map $\cM\to S$ which is a contradiction. Otherwise, $\pi_S(\bP^1_E)$ defines a rational curve on $S$ moving with $E$, contradicting the rigidity of rational curves on $S$.

Now let us show that the intersection is nonempty (then it will be automatically $2$-dimensional). For the general $E\in\cM$, consider the degree 2 map $D_E\to\bP^1_E$ where
\[
D_E:=\{(\xi,p)\in \bP^1_E\times S\;|\;p\in\Supp(\xi)\}
\] 
We claim that this degree 2 map must have a branch point. If $D_E$ is irreducible, this follows from Hurwitz formula. Otherwise, $D_E$ consists of two rational components; thus the image of the natural map $D_E\to S$ produces two rational curves on $S$. Since rational curves on $S$ are rigid, it follows that this image (moving $E$) is a fixed curve $C$ with (at most) two rational components. But the locus of length 2 subschemes supported at $C$ is 2-dimensional, hence $R_2$ can't be contained in this locus (recall that $R_2$ is 3-dimensional), which gives the contradiction.

Now consider a component $W_S\subset R_2\cap S^{[2]}_{nred}$ such that the projection $\pi_\cM:W_S\to \cM$ is surjective. This exists, since we just proved that $D_E$ is irreducible for $E\in\cM$ general (hence $\bP^1_E$ contains a nonreduced subscheme). The general fiber of the map $\pi_S:W_S\to S$ consists of nonreduced subschemes supported at a general point $p$, hence it is contained in a rational curve. We deduce that if the general fiber of $\pi_S$ is positive-dimensional, then $W_S$ is covered by rational curves (the fibers of $\pi_S$); hence $\cM$ is covered by rational curves (recall that $W_S\to\cM$ is dominant), contradiction. It follows that $W_S\to S$ is surjective, as desired. 
\end{proof}

Now we give the new construction of rational maps of degree $4$:

\begin{lem}
 There is a  morphism $W_S\to W^2_4(S,L)$ generically finite onto its image.
 \label{deg44}
\end{lem}

\begin{proof}
For any (nonreduced) $\xi\in W_S$ we have an associated $E$ with $h^1(E\otimes \cI_\xi)=1$. We may consider the extension
\[
0 \longrightarrow E^\vee \longrightarrow F \longrightarrow \cI_\xi \longrightarrow 0
\]
and the subsheaf $E_1:=\mathrm{Im}(F^\vee \to E)\subset E$. By numerical reasons we get $\mathrm{length}(E/E_1)=c_2(F^\vee)-c_2(E)=2$, and also $E\otimes \cI_\xi \subset E_1$; indeed, by  applying the functor $\sHom(-,\cO_S)$ to the sequence above we get that $E_1$ is the kernel of $E\to \sExt^1(\cI_\xi,\cO_S)=\cO_\xi$.
Moreover, if $p=\Supp(\xi)$, then $E_1\otimes k(p)$ is $3$-dimensional since $E_1$ is a quotient of a rank 3 locally free sheaf. It is then easy to deduce that $(E_1)_p=e_1\cdot\cO_{S,p}+e_2 \cdot\cI_\xi$ for some $e_1,e_2$ local generators of $E_p$. 

Note that $h^1(E_1)=1$ (since $h^1(F^\vee)=h^2(F^\vee)=0$). 
It follows that the subsheaf $\tilde E\subset E_1$ defined as $\tilde E_p=e_1\cdot m_p+e_2\cdot m_p^2$ around $p$ (and $\tilde E_q=E_q$ for $q\neq p$) satisfies $h^1(\tilde E)\ge 1$ (since $\tilde E\subset E_1$). Actually, a quick numerical check gives $h^1(\tilde E)= 1$ and $\chi(\tilde E)=2$, hence $h^0(\tilde E)=3$. The subspace $H^0(\tilde E)\subset H^0(E)$ induces a map of degree $\leq 4$ by \autoref{compdegree} (all sections are vanishing with order $2$ at $p$).

The corresponding morphism $W_S\to W^2_4(S,L)$ has image of dimension $2$, since $\pi_\cM:W_S \to \cM$ is surjective and equals the composition
\[
W_S \longrightarrow W^2_4(S,L)_6 \longrightarrow \cM,
\]
where $W^2_4(S,L)_6$ parametrizes maps of degree $\le 4$ having kernel bundle with $c_2=6$ (this is the notation in the preliminaries), and the right arrow is the forgetful map  $(E,V^\vee)\mapsto E$.
\end{proof}

Now we show that the image of $W_S\to W^2_4(S,L)$ is not contained in $R_2\subset W^2_4(S,L)$, so that the above construction really yields new maps of degree $4$:

\begin{prop}
    The locus $\mathrm{Im}(W_S)\subset W^2_4(S,L)$ is not contained in $R_2\subset W^2_4(S,L)$.
\end{prop}
\begin{proof}
Recall that a rational map $S\dashrightarrow \mathbb P^2$ induced by the primitive linear system is completely determined by the couple $(E,V^\vee)$. We have $(E,V^\vee)\in R_2$ if and only if $V^\vee =H^0(E\otimes \cI_\zeta)$ for some $\zeta\in S^{[2]}$. Assume that $H^0(\tilde E)=H^0(E\otimes \cI_\zeta)$, where (following the notation of \autoref{deg44}) $\tilde E$ is constructed from a nonreduced $\xi\in W_S$.

If $\zeta$ is reduced, then all sections of $H^0(E\otimes \cI_\zeta)$ vanish with order $2$ at the point where $\xi$ is supported. Hence $\bigcap _{s\in V^\vee}Z_{cycle}(s)$ is of degree $\geq3$, which by \autoref{compdegree} implies that $(E,V^\vee)\in W^2_3(S,L)$, in contradiction to \autoref{ccc}. 

Now suppose that $\zeta$ is nonreduced. We get that $\zeta$ is supported at the point $p$ where $\tilde E$ is not locally free. This implies that all the sections in $V^\vee$ are locally contained in the subsheaf $e_1\cdot \cI_\zeta +e_2\cdot m_p^2$ of $E$; in particular, we get $h^0(E\otimes \cI_p^2)\ge 2$. In order to finish the proof, it suffices to prove the following claim: for a general point $p\in S$ there is no vector bundle $E_p\in \cM$ such that $h^0(E_p\otimes \cI_p^2)=h^1(E_p\otimes \cI_p^2)=2$.

Observe that such a bundle $E_p\in\cM$ is a singular point of the curve
\[C_{\cI_p^2}=\{E\in \cM \;|\; h^1(E\otimes \cI_p^2)\ge 1 \},
\]
which has arithmetic genus $3$ (recall \autoref{coh913}).
Therefore, if the claim is false one can consider the rational map $S\dashrightarrow \cM$ sending a general $p\in S$ to the unique singular point of $C_{\cI_p^2}$ (note that if $C_{\cI_p^2}$ has at least two singular points for every $p\in S$, then $\{C_{\cI_p^2}\}_{p\in S}$ defines a $2$-dimensional family of curves on $\cM$ with geometric genus $\leq 1$, contradiction).
This rational map actually extends to a morphism $S\to \cM$; indeed, the finite birational morphism
\[
\{(p,E)\in S\times\cM\;|\;h^0(E\otimes\cI_p^2)\geq 2\}\longrightarrow S
\]
must be an isomorphism, by normality of $S$. Now the morphism $S\to \cM$ is étale by Hurwitz formula, which contradicts the fact that $\cM$ is simply connected. This finishes the proof.
\end{proof}

\section{Genus 11}\label{sec:genus11}

Let $(S,L)$ be a polarized $K3$ surface of genus 11 with $\Pic(S)\cong\bZ\cdot L$, and let $\cM$ denote the (fine) moduli space $\cM(2,L,5)$. Any $E\in\cM$ has $h^0(E)=7$ and $c_2(E)=7$. For every $\xi\in S^{[2]}$, considering $H^0(E\otimes\cI_\xi)\subset H^0(E)$ one can construct maps $S\dasharrow\bP^2$ of degree $\leq5$. In particular, $W^2_5(S,L)$ has (a birational copy of) $S^{[2]}\times\cM$ as a component. Its singular locus determines the component of $W^2_4(S,L)$ that we will describe, see \autoref{singu}.

\subsection{Constructing maps of degree 4}
Our first goal is to construct rational maps of degree 4 (and to determine the component of $W^2_4(S,L)$ obtained from this construction). We can do this thanks to:

\begin{prop}\label{cohgenus11}
    For every $E\in\cM$, there exists a unique $\xi\in S^{[3]}$ such that $h^1(E\otimes\cI_\xi)\geq2$. In particular, $W^2_4(S,L)$ admits an irreducible component isomorphic to $\cM$.
\end{prop}

The arguments for proving this proposition rely on Bridgeland stability. Along the proof, we will determine several explicit isomorphisms between moduli spaces of Gieseker stable sheaves on $S$ and punctual Hilbert schemes on $\cM$. In particular, we will construct an isomorphism
\[
\cM^{[3]}\overset{\cong}{\longrightarrow} \cM(1,0,-2)=S^{[3]}
\]
Under this isomorphism, the unique $\xi\in S^{[3]}$ such that $h^1(E\otimes\cI_\xi)\geq2$ corresponds to the unique non-curvilinear element of $\cM^{[3]}$ which is supported on $E$.

Consider $\xi\in S^{[3]}$ and $E\in\cM$. The numerical wall $W$ defined by $E^\vee[1]$ and $\cI_\xi$ has endpoints $\beta=\frac{-1}{2}$ and $\beta=\frac{-2}{5}$. Since $\cI_\xi$ (resp. $E^\vee[1]$) is $\sigma_{\alpha,\frac{-1}{2}}$-stable (resp. $\sigma_{\alpha,\frac{-2}{5}}$-stable) for every $\alpha>0$ by \autoref{minimalrank}, it follows that both $\cI_\xi$ and $E^\vee[1]$ are semistable along $W$.

Furthermore, there is a hole along $W$ created by the spherical bundle $G$ with $v(G)=(7,-3L,13)$ (so that $W=W(\cI_\xi, E^\vee[1])=W(G,E^\vee[1])$). Finally, we also note that both $E^\vee[1]$ and $G$ are $\sigma_{\alpha,\frac{-4}{9}}$-stable for every $\alpha>0$; we will write $\sigma_0:=W\cap\{\beta=\frac{-4}{9}\}$.

As usual, any nonzero morphism $\cI_\xi\to E^\vee[1]$ is a surjection in $\Coh^\beta(S)$ ($\frac{-1}{2}\leq\beta<\frac{-1}{3}$), defining a destabilization of $\cI_\xi$ along $W$. In order to understand the possible number of morphisms $\cI_\xi\to E^\vee[1]$, we also need to control the possible HN filtrations of $\cI_\xi$ below $W$.

As in genus 7, we achieve this by establishing several isomorphisms of (Hilbert schemes of) moduli spaces. For instance, the same arguments as in \autoref{isomoduli} prove:

\begin{lem}
    There is an isomorphism $\cM\overset{\cong}{\longrightarrow}\cM(5,-2L,8)$ sending a bundle $E\in\cM$ to the (unique up to isomorphism) sheaf $P$ fitting in a nontrivial extension
    \[
    0\to G\to P\to E^\vee[1]\to 0.
    \]
    The inverse is obtained by sending $P\in\cM(5,-2L,8)$ to (the shift of the derived dual of) its last HN factor for stability conditions just below $\sigma_0$.
\end{lem}

\begin{rem}
    By rotating the distinguished triangle, $0\to P\to E^\vee[1]\to G[1]\to 0$ defines a destabilization of $E^\vee[1]$ along $W\cap\{\beta>\frac{-3}{7}\}$.
\end{rem}

With similar arguments to those of \autoref{isohilb2}, one obtains:

\begin{lem}
    There is an isomorphism $\cM^{[2]}\overset{\cong}{\longrightarrow}\cM(3,-L,3)$.
\end{lem}

Finally, one more step in this construction yields our desired isomorphism of Hilbert schemes:

\begin{lem}\label{isohilb3}
    There is an isomorphism $\gamma:\cM^{[3]}\overset{\cong}{\longrightarrow}\cM(1,0,-2)=S^{[3]}$.
\end{lem}
\begin{proof}
    The proof is similar to that of \autoref{isohilb2}. We define the morphism $\gamma$ as follows.

    Given a reduced element of $\cM^{[3]}$, supported at $E_1,E_2,E_3\in\cM$, we consider the object $M$ fitting into an extension
\[
0\to G\to M\to E_1^\vee[1]\oplus E_2^\vee[1]\oplus E_3^\vee[1]\to 0
\]
such that the induced elements in $\Ext^1(E_i^\vee[1],G)$ are nonzero. Such an $M$ is unique up to isomorphism (due to the action of automorphisms of $E_1^\vee[1]\oplus E_2^\vee[1]\oplus E_3^\vee[1]$), and of the form $M=\cI_\xi$ for some $\xi\in S^{[3]}$ (indeed $M$ is semistable above $\sigma_0$, by arguing as in \autoref{isohilb2}).

If our element of $\cM^{[3]}$ consists of a tangent direction $T\in\bP(\Ext^1(E_1,E_1))$ at $E_1$ and a reduced point $E_2$, we consider $\widetilde{T}:=T^\vee[1]\in\bP(\Ext^1(E_1^\vee[1],E_1^\vee[1]))$ and take the object $M$ fitting into
\[
0\to G\to M\to \widetilde{T}\oplus E_2^\vee[1]\to 0,
\]
such that the induced elements in $\Ext^1(E_1^\vee[1],G)$ and $\Ext^1(E_2^\vee[1],G)$ (the former via the map $\Ext^1(\widetilde{T}\oplus E_2^\vee[1],G)\to\Ext^1(\widetilde{T},G)\to\Ext^1(E_1^\vee[1],G)$) are nonzero. Again, $M$ is unique up to isomorphism and of the form $M=\cI_\xi$ for some $\xi\in S^{[3]}$.

Finally, given a tangent direction $T\in\bP(\Ext^1(E^\vee[1],E^\vee[1]))$ at $E\in\cM$, it is well known that all elements of $\cM^{[3]}$ supported at $E$ and containing $T$ are parametrized by $\bP(\Ext^1(T,E))$. Given such a $Q\in\bP(\Ext^1(T,E))$, consider the corresponding $\widetilde{Q}:=Q^\vee[1]\in\bP(\Ext^1(E^\vee[1],T^\vee[1]))$ and take the object $M$ sitting in an extension
\[
0\to G\to M\to \widetilde{Q}\to 0,
\]
such that the induced element in $\Ext^1(E^\vee[1],G)$ (via the map $\Ext^1(\widetilde{Q},G)\to\Ext^1(T^\vee[1],G)\to\Ext^1(E^\vee[1],G)$) is nonzero. Again $M$ is unique and of the form $M=\cI_\xi$, for some $\xi\in S^{[3]}$.

The morphism $\gamma$ is injective on closed points and hence, arguing as in the proof of \autoref{isohilb2}, we obtain that $\gamma$ is an isomorphism.
\end{proof}

Now we have all the ingredients to prove \autoref{cohgenus11}:

\begin{proof}[Proof of \autoref{cohgenus11}]
Given $\xi\in S^{[3]}$ and $E'\in\cM$, we want to determine the number $\hom(\cI_\xi, E'^\vee[1])=h^1(E'\otimes\cI_\xi)$. To this end, we consider the HN filtration
\[
0\to G\to \cI_\xi\to F\to 0
\]
of $\cI_\xi$ for stability conditions just below $\sigma_0$, established in \autoref{isohilb3}. The object $F$ is nothing but $F=\Phi_\cE(\cO_{\gamma^{-1}(\xi)})^\vee[1]$, where $\Phi_\cE:\Db(\cM)\to \Db(S)$ is the Fourier-Mukai equivalence with kernel the universal bundle $\cE$ on $S\times\cM$.

Due to the vanishing $\Hom(G,E'^\vee[1])=0$, we have equalities
\[
\begin{aligned}
&\Hom_{\Db(S)}(\cI_\xi, E'^\vee[1])=\Hom_{\Db(S)}(F, E'^\vee[1])=\\
&=\Hom_{\Db(S)}(\Phi_\cE(\cO_{\gamma^{-1}(\xi)})^\vee[1], \Phi_\cE(k(E'))^\vee[1])=\Hom_{\Db(\cM)}(k(E'),\cO_{\gamma^{-1}(\xi)}).
\end{aligned}
\]

If $E'\notin \Supp \gamma^{-1}(\xi)$, it immediately follows that $\hom(\cI_\xi, E'^\vee[1])=0$. 

On the other hand, if $E'\in \Supp \gamma^{-1}(\xi)$, then the equality 
\[
\hom(k(E'),\cO_{\gamma^{-1}(\xi)})=\ext^1(k(E'),\cO_{\gamma^{-1}(\xi)})-1=\ext^1(\cO_{\gamma^{-1}(\xi)},k(E'))-1
\]
holds; indeed, it is a consequence of the equalities $\ext^2(k(E'),\cO_{\gamma^{-1}(\xi)})=\hom(\cO_{\gamma^{-1}(\xi)},k(E'))=1$ and $\chi(k(E'),\cO_{\gamma^{-1}(\xi)})=0$.

Therefore, by applying \autoref{numbergen} we obtain
    \[
    \hom(\cI_\xi, E'^\vee[1])=\left\{
    \begin{array}{c l}
     1 & \text{if $\gamma^{-1}(\xi)$ is curvilinear at $E'$}\\
     2 & \text{if $\gamma^{-1}(\xi)$ is not curvilinear at $E'$}\\
    \end{array}
    \right.
    \]

In particular, for a given $E\in\cM$, there exists a unique $\xi\in S^{[3]}$ such that $h^1(E\otimes\cI_\xi)=2$; it is nothing but the image under $\gamma$ of the unique non-curvilinear length 3 subscheme supported at $E\in\cM$. This concludes the proof of the proposition.
\end{proof}

\subsection{Discarding maps of degree 3 (I)} 

In view of \autoref{minimalc2}, every degree 3 map $S\dasharrow \bP^2$ arises from $V^\vee\in\Gr(3,H^0(E))$, where $E\in\cM$ and $\deg\left(\bigcap_{s\in  V^\vee}Z_{cycle}(s)\right)=4$. In this subsection we discard the case where the schematic intersection $\bigcap_{s\in  V^\vee}Z(s)$ is of length 4, by means of the following:

\begin{prop}\label{coh2genus11}
    For all $E\in\cM$ and $\xi\in S^{[4]}$, the inequality $h^0(E\otimes\cI_\xi)\leq 2$ holds.
\end{prop}

In order to prove this result, we perform a Bridgeland analysis which has many symmetries with that of the previous subsection.
In this case, the numerical wall $\overline{W}$ defined by $E^\vee$ and $\cI_\xi$ has endpoints $\beta=\frac{-3}{5}$ and $\beta=\frac{-1}{2}$. Any nonzero morphism $E^\vee\to\cI_\xi$ is an injection in $\Coh^\beta(S)$ ($\frac{-3}{5}\leq\beta<\frac{-1}{2}$), which destabilizes $\cI_\xi$ along $\overline{W}$.

The wall $\overline{W}$ has a hole created by the spherical bundle $\overline{G}$ with $v(\overline{G})=(7,-4L,23)$. Moreover, $E^\vee,\overline{G}[1]\in\Coh^{\frac{-5}{9}}(S)$ are $\sigma_{\alpha,\frac{-5}{9}}$-stable for every $\alpha>0$; let us write $\overline{\sigma_0}:=\overline{W}\cap\{\beta=\frac{-5}{9}\}$.

Again, by exploiting the $\overline{\sigma_0}$-stability of $E^\vee$ and $\overline{G}[1]$ and the fact that $\ext^1(\overline{G}[1],E^\vee)=1$, one can construct several isomorphisms of (Hilbert schemes of) moduli spaces:

\begin{lem} Let $\sigma'$ be a stability condition above $\overline{W}$, close enough to $\overline{\sigma_0}$. Then there are isomorphisms:
    \begin{enumerate}[{\rm (1)}]
    \item $\cM(2,-L,5)\overset{\cong}{\longrightarrow}\cM_{\sigma'}(-5,3L,-18)$.
    \item $\cM(2,-L,5)^{[2]}\overset{\cong}{\longrightarrow}\cM_{\sigma'}(-3,2L,-13)$.
    \item $\cM(2,-L,5)^{[3]}\overset{\cong}{\longrightarrow}\cM_{\sigma'}(-1,L,-8)$.
    \item  $\overline{\gamma}:\cM(2,-L,5)^{[4]}\overset{\cong}{\longrightarrow}\cM_{\sigma'}(1,0,-3)=S^{[4]}$.
\end{enumerate}
\end{lem}
\begin{proof}
    The arguments are similar to those of the previous subsection. For instance, in (1) one sends $M\in\cM(2,-L,5)$ to the isomorphism class of the object $P$ sitting in a non-trivial extension $0\to M\to P\to \overline{G}[1]\to 0$.
\end{proof}

\begin{rem}
    Let $(S,L)$ be a polarized $K3$ surface of genus $g$, $\Pic(S)=\bZ\cdot L$. The isomorphisms of Hilbert schemes $S^{[2]}\cong\cM(2,L,3)^{[2]}$ (in genus 7) and $S^{[3]}\cong\cM(2,L,5)^{[3]}$, $S^{[4]}\cong\cM(2,L,5)^{[4]}$ (in genus 11) can be generalized, with similar arguments, to the following statement: if $g=4n-1$, then there are isomorphisms $S^{[n]}\cong\cM(2,L,2n-1)^{[n]}$, $S^{[n+1]}\cong\cM(2,L,2n-1)^{[n+1]}$.

    This can be thought of as a reinterpretation, in terms of Bridgeland stability, of an isomorphism of moduli spaces of sheaves on $S$ and a Fourier-Mukai partner, which is due to Yoshioka and holds in broader generality (\cite[Theorem 7.6]{yos}).
\end{rem}

Now we can prove \autoref{coh2genus11} with similar arguments to those of \autoref{cohgenus11}:

\begin{proof}[Proof of \autoref{coh2genus11}]
Given $\xi\in S^{[4]}$ and $E'\in\cM$, we want to determine the number $\hom(E'^\vee,\cI_\xi)=h^0(E'\otimes\cI_\xi)$. We will use the isomorphism $\overline{\gamma}$ of the previous lemma, together with the fact that $E'^\vee\in\cM(2,-L,5)$. Consider the HN filtration
\[
0\to \overline{F}\to \cI_\xi\to \overline{G}[1]\to 0
\]
of $\cI_\xi$ for stability conditions just below $\overline{\sigma_0}$; the factor $\overline{F}$ is determined by $\overline{\gamma}^{-1}(\xi)$. Indeed, one has $\overline{F}=\Phi_{\overline{\cE}}(\cO_{\overline{\gamma}^{-1}(\xi)})$, where $\overline{\cE}$ denotes the universal bundle on $S\times \cM(2,-L,5)$ and $\Phi_{\overline{\cE}}:\Db(\cM(2,-L,5))\to \Db(S)$ is the corresponding Fourier-Mukai equivalence.

Due to the vanishing $\Hom(E'^\vee,\overline{G}[1])=0$, we have 
\[
\begin{aligned}
&\Hom_{\Db(S)}(E'^\vee,\cI_\xi)=\Hom_{\Db(S)}(E'^\vee, \overline{F})=\Hom_{\Db(S)}(\Phi_{\overline{\cE}}(k(E'^\vee), \Phi_{\overline{\cE}}(\cO_{\overline{\gamma}^{-1}(\xi)}))=\\
&=\Hom_{\Db(\cM(2,-L,5))}(k(E'^\vee),\cO_{\overline{\gamma}^{-1}(\xi)}).
\end{aligned}
\]

If $E'^\vee\notin \Supp \overline{\gamma}^{-1}(\xi)$, the vanishing $\hom(E'^\vee,\cI_\xi)=0$ follows. 

On the other hand, if $E'^\vee\in \Supp \overline{\gamma}^{-1}(\xi)$, we have the equality 
\[
\hom(k(E'^\vee),\cO_{\overline{\gamma}^{-1}(\xi)})=\ext^1(k(E'^\vee),\cO_{\overline{\gamma}^{-1}(\xi)})-1=\ext^1(\cO_{\overline{\gamma}^{-1}(\xi)},k(E'^\vee))-1
\]
(following from $\ext^2(k(E'^\vee),\cO_{\overline{\gamma}^{-1}(\xi)})=\hom(\cO_{\overline{\gamma}^{-1}(\xi)},k(E'^\vee))=1$ and $\chi(k(E'^\vee),\cO_{\overline{\gamma}^{-1}(\xi)})=0$).
Therefore, by applying \autoref{numbergen} we deduce that
    \[
    \hom(E'^\vee, \cI_\xi)=\left\{
    \begin{array}{c l}
     1 & \text{if $\overline{\gamma}^{-1}(\xi)$ is of complete intersection at $E'^\vee$}\\
     2 & \text{if $\overline{\gamma}^{-1}(\xi)$ is not of complete intersection at $E'^\vee$}\\
    \end{array}
    \right.
    \]
which finishes the proof.
\end{proof}

Furthermore, the analysis performed in the proof of \autoref{coh2genus11} yields the following:

\begin{cor}\label{jumplocus}
Let $E\in \cM$. Then:
\begin{enumerate}[{\rm (1)}]
    \item The jump locus $\{\xi \in S^{[4]}\;|\;h^0(E\otimes \cI_\xi)\ge 2\}$ is isomorphic to $\mathrm{Bl}_{E}(\cM)$.
    \item Under this isomorphism, every $\xi\in S^{[4]}$ in the exceptional divisor satisfies $h^0(E'\otimes \cI_\xi)=0$ for all $E'\in\cM\setminus\{E\}$.
\end{enumerate}
\end{cor}

\subsection{Discarding maps of degree 3 (II)}
In this subsection we complete the proof of the emptiness of $W^2_3(S,L)$. Assume given a map $\phi_V:S\dasharrow \bP^2$ of degree 3, defined by $V^\vee\in\Gr(3,H^0(E))$ for $E\in\cM$. The base ideal $\cI$ of the linear system $V\subset H^0(L)$ has colength $L^2-c_2(E)=13$, and $\deg\left(\bigcap_{s\in  V^\vee}Z_{cycle}(s)\right)=4$; according to \autoref{coh2genus11}, the schematic intersection $\xi:=\bigcap_{s\in  V^\vee}Z(s)$ must be of length $\leq 3$. 

\begin{lem}
    Under these assumptions, the schematic intersection $\xi$ has length 3.
\end{lem}
\begin{proof}
    If $\xi$ consists of a reduced point $p\in S$, then by \autoref{ideallemma} there is a local trivialization $E_p\cong e_1\cdot \cO_{S,p}\oplus e_2\cdot\cO_{S,p}$ of $E$ near $p$, such that $s_p\in m_pe_1+m_p^4 e_2$ for every $s\in V^\vee$. This implies $\cI\subset \cI_p^5$ and hence $\colength(\cI)\geq\colength(\cI_p^5)=15$, which is a contradiction.

    Assume now that $\xi$ is of length 2, supported at two distinct points $p$ and $q$. If $\bigcap_{s\in  V^\vee}Z_{cycle}(s)=2p+2q$, then by \autoref{ideallemma} one has $s_p\in m_pe_1+m_p^2e_2$ for all $s\in V^\vee$, in an appropriate trivialization of $E$ near $p$; therefore $\cI\subset \cI_p^3$, and since $\cI$ has (at most) three local generators at $p$, it follows that $\colength(\cI,p)\geq 7$. Similarly we have $\colength(\cI,q)\geq 7$, hence $\colength(\cI)\geq 14$ which is impossible.

    On the other hand, if $\bigcap_{s\in  V^\vee}Z_{cycle}(s)=3p+q$, a similar use of \autoref{ideallemma} yields $\cI\subset\cI_p^4\cap\cI_q^2$; since $\cI$ has three generators at $p$, we have $\colength(\cI,p)\geq 12$ and therefore $\colength(\cI)\geq\colength(\cI,p)+\colength(\cI,q)\geq 15$, contradiction.

    It only remains to discard the situation where $\xi$ is nonreduced of length 2, supported at a point $p$. In this case, by \autoref{ideallemma} we can find (in appropriate local coordinates) generators $e_1,e_2$ of $E$ around $p$, such that $s_p\in(x,y^2)e_1+(x^2,xy^2,y^4)e_2$ for every $s\in V^\vee$. Under the isomorphism $V\cong \bigwedge^2V^\vee$, this implies that $\cI\subset (x^3,x^2y^2,xy^4,y^6)$ locally at $p$; since $\cI$ has at most three generators at $p$, it follows that $\colength(\cI,p)\geq 14$, which gives a contradiction.
\end{proof}

Apart from having a schematic intersection $\xi=\bigcap_{s\in  V^\vee}Z(s)$ of length 3, we can reduce to a more concrete situation:

\begin{lem}
    Let $p\in S$ be the point where $\mathrm{cycle}(\xi)$ and $\bigcap_{s\in  V^\vee}Z_{cycle}(s)$ differ. Then $\xi$ is reduced at $p$.
\end{lem}
\begin{proof}
    Let us first show that $\xi$ is curvilinear. Assume for the sake of a contradiction that $\xi$ is not curvilinear, namely $\cI_\xi=\cI_p^2$. Since $h^0(E\otimes\cI_\xi)\geq 3$, we have $h^1(E\otimes\cI_\xi)\geq 2$ and we can consider a stable sheaf $F$ arising as a double extension
    \[
    0\longrightarrow\cO_S^{\oplus 2}\longrightarrow F\longrightarrow E\otimes\cI_p^2\longrightarrow 0.
    \]

    Note that $v(F)=(4,L,1)$. Since $E\otimes\cI_p^2\otimes k(p)$ is a 6-dimensional ($E\otimes\cI_p^2$ has 6 local generators at $p$), it follows that the minimal number of generators of $F$ at $p$ is at least 6. This implies that the reflexive hull $F^{**}$ of $F$ has Mukai vector $v(F^{**})=(4,L,k)$ with $k\geq 3$, hence $v(F^{**})^2\leq -4$ which contradicts the stability of $F^{**}$.

    Now we only need to rule out the scenario where $\xi$ is of the form $(x,y^m)$ locally at $p$, for $m\in\{2,3\}$. In that case, by \autoref{ideallemma} we have $s_p\in(x,y^m)\cdot e_1+(x^2,xy,y^{m+1})\cdot e_2$ for all $s\in V^\vee$, in an appropriate trivialization $E_p\cong e_1\cdot \cO_{S,p}\oplus e_2\cdot \cO_{S,p}$. In view of the isomorphism $V\cong\bigwedge^2V^\vee$, it follows that the base ideal $\cI$ satisfies $\cI\otimes\hat \cO_{S,p}\subset (x^3,x^2y, xy^{m+1},y^{2m+1})$ locally at $p$; furthermore, recall that $\cI$ has at most three local generators at $p$.

    For $m=3$, this implies that $\colength(\cI,p)\geq 14$, which is a contradiction.

    For $m=2$, this implies that $\colength(\cI,p)\geq 11$. If $q\neq p$ is the other point where $\xi$ is supported, then $\colength(\cI,q)\geq 3$; hence $\colength(\cI)\geq 14$, and again we obtain a contradiction.
\end{proof}

The last ingredient we need to derive a contradiction from this situation is the following lemma:

\begin{lem}\label{minimalclass}
Let $\xi\in S^{[3]}$ be a length 3 subscheme, which is reduced at  $p\in S$. Then the locus
\[
C_{p,\xi}:=\{E'\in\cM\;|\;\exists s\in H^0(E'\otimes\cI_\xi)\setminus\{0\}\text{ such that $Z(s)$ has length $\geq 2$ at $p$}\}
\]
describes a curve in the linear system $|\hat{L}|$, where $\Pic(\cM)= \bZ\cdot \hat{L}$.
\end{lem}
\begin{proof}
Let us denote by $\zeta$ the (length 2) residual subscheme of $p$ in $\xi$. Let $\cE$ be a universal bundle on $S\times\cM$, defining a Fourier-Mukai equivalence $\widetilde{\Phi_\cE}:=\Phi_\cE^{S\to \cM}:\Db(S)\to \Db(\cM)$.

The first observation is that, for an appropriate choice of $\cE$ (obtained by twisting $\cE$ with the pullback of a line bundle on $\cM$), we have $R^0\widetilde{\Phi_\cE}(\cI_\xi)\cong \cO_\cM$. This can be deduced from the following facts:
\begin{itemize}
    \item $v(\widetilde{\Phi_\cE}(\cI_\xi))=(1,0,-2)$ (see e.g. \cite[Lemma~7.2]{yos} or \cite[Proposition~15]{Maciocia}).

    \item $R^1\widetilde{\Phi_\cE}(\cI_\xi)=\sExt^2(\cO_{\gamma^{-1}(\xi)},\cO_\cM)$, where $\gamma:\cM^{[3]}\to S^{[3]}$ is the isomorphism of \autoref{isohilb3}. Indeed, to prove this equality one simply applies the Fourier-Mukai functor $\widetilde{\Phi_\cE}$ to the distinguished triangle
    \[
    0\to G\to \cI_\xi\to \Phi^{\cM\to S}_\cE(\cO_{\gamma^{-1}(\xi)})^\vee[1]\to 0
    \]
    (see the proof of \autoref{cohgenus11}), and uses that $\Phi^{\cM\to S}_\cE(\cdot^\vee)^\vee$ is the quasi-inverse of $\widetilde{\Phi_\cE}$.
\end{itemize}

The unique (up to constant) global section of $R^0\widetilde{\Phi_\cE}(\cI_\xi)\cong \cO_\cM$ restricts, at a point $E'\in\cM\setminus\gamma^{-1}(\xi)$, to the unique section of $H^0(S,E'\otimes\cI_\xi)$. We will denote by $\widetilde{s}\in H^0(S\times\cM,\pi_1^*(\cI_\xi)\otimes\cE)$ the corresponding section, where $\pi_1:S\times\cM\to S$ is the first projection.

Locally at a point $E'\in\cM$, the locus $C_{p,\xi}$ is given as follows. Fix $x,y$ local analytic coordinates of $S$ around $p$, $z,w$ local coordinates of $\cM$ around $E'$ and a local trivialization of $\cE$ around $(p,E')\in S\times\cM$. Then $\tilde s$ has a local expression
\[
(a_1(z,w)\cdot x+b_1(z,w)\cdot y+O(x,y)^2,a_2(z,w)\cdot x+b_2(z,w)\cdot y+O(x,y)^2)
\]
and $C_{p,\xi}$ is locally given by the equation $a_1b_2-a_2b_1=0$.

To determine the class of the curve $C_{p,\xi}$, we pick a basis $w_1,w_2$ of the tangent space $T_pS$ at $p$. For $i=1,2$ the section $\tilde s$ gives rise to a section
\[
\tilde{s_i}\in H^0(S\times\cM,\pi_1^*\left(\cI_{\xi}/(\cI_\zeta\otimes\cI_{[w_i]})\right)\otimes\cE)=H^0(\cM,R^0\widetilde{\Phi_\cE}\left(\cI_{\xi}/(\cI_\zeta\otimes\cI_{[w_i]})\right)).
\]

Now observe that $\cI_{\xi}/(\cI_\zeta\otimes\cI_{[w_i]})$ is isomorphic (as an $\cO_S$-module) to the skyscraper sheaf $k(p)$; this isomorphism is canonical (up to constant), defined via $w_i$.
We obtain an isomorphism
\[
R^0\widetilde{\Phi_\cE}\left(\cI_{\xi}/(\cI_\zeta\otimes\cI_{[w_i]})\right)\cong R^0\widetilde{\Phi_\cE}\left(k(p)\right)=\widetilde{\Phi_\cE}(k(p))=\cE_p
\]
given by $w_i$, where $\cE_p\cong\cE|_{\{p\}\times\cM}$ is the vector bundle on $\cM$ defined by the point $p\in S$.

The curve $C_{p,\xi}$ is then defined as the zero locus of $\tilde{s_1}\wedge\tilde{s_2}\in H^0(\bigwedge ^2 \cE_p)$. Since $v(\cE_p)=(2,\hat{L},5)$  (see again \cite[Proposition~15]{Maciocia})), we have $\bigwedge ^2 \cE_p=\hat{L}$ and the assertion follows.
\end{proof}

Now we can conclude the proof of the nonexistence of degree 3 maps. For this last argument, the assumption $\Pic(S)=\bZ\cdot L$ is not enough; we will also assume that $(S,L)$ is general among $K3$ surfaces of Picard rank 1.

\begin{prop}
If $(S,L)$ is general among $K3$ surfaces of Picard rank 1, then there exists no $V^\vee \in \mathrm{Gr}(3,H^0(E))$ such that $\bigcap_{s\in V^\vee} Z_{cycle}(s)=2p+q+r$ (with $p\neq q,r$) and $\bigcap_{s\in V^\vee}Z(s)$ is of length 3 and reduced at $p$. 
\end{prop}
\begin{proof}
We argue by contradiction.
As in the proof of \autoref{minimalclass}, we denote by $\zeta$ the (length 2) residual subscheme of $p$ in $\bigcap_{s\in V^\vee}Z(s)$; in the notations of the statement, we have $\zeta=q+r$ with $q,r$ two (possibly infinitely near) points.

First notice that for any length 2 subscheme $\eta\in \bP(T_pS)$ nonreduced at $p$, the inequality $h^0(E\otimes \cI_{\eta\cup\zeta})\ge 2$ holds; this follows from \autoref{ideallemma}. 
Hence the subschemes $\{\eta\cup\zeta\;|\;\eta\in\bP(T_pS)\}$ describe a smooth rational curve in the jump locus $\mathrm{Bl}_{E}(\cM)\subset S^{[4]}$ of \autoref{jumplocus}. We claim that this curve is the exceptional divisor in $\mathrm{Bl}_{E}(\cM)$; in particular, by \autoref{jumplocus}.(2) this will imply that $h^0(\tilde E \otimes \cI_{\eta\cup\zeta})=0$ for any $\eta\in \bP(T_pS)$ and $\tilde E \in\cM\setminus\{ E\}$.

In order to prove the claim, assume that the curve is not the exceptional divisor. Then its image $R$ in $\cM$ is a rational curve, describing a component of the curve
\[
C_{p,\xi}:=\{E'\in\cM\;|\;\exists s\in H^0(E'\otimes\cI_\xi)\setminus\{0\}\text{ such that $Z(s)$ has length $\geq 2$ at $p$}\}.
\]
The curve $C_{p,\xi}$ lies in the primitive linear system $|\hat{L}|$ of $\cM$ by \autoref{minimalclass}, hence $R=C_{p,\xi}\in |\hat{L}|$ and $p_a(R)=11$. This yields a contradiction, since on the one hand $R$ can only be singular at $E\in\cM$, whereas on the other hand every rational curve in $|\hat{L}|$ is nodal under our generality assumptions on $(S,L)$ (hence on $(\cM,\hat{L})$) by Chen's result \cite{chen}.

Now we start the analysis.
Under our hypotheses we get $V=\bigwedge ^2 V^\vee\subset H^0(L\otimes \cI_p^3\otimes \cI_\zeta^2)$. Given $z,w$ local coordinates at $p$, consider for every $\lambda\in\bC$ the colength 4 ideal $\mathcal{J}_\lambda:=(z(z+\lambda w),w(z+\lambda w),w^3)\supset \cI_p^3$. Since $V\subset H^0(L \otimes \mathcal J_\lambda \otimes \cI_\zeta^2)$, we get $h^1(L \otimes \mathcal J_\lambda \otimes \cI_\zeta^2)\ge 1$.  

Given $\lambda\in\bC$, we consider any nontrivial extension $\tilde E_\lambda$ of $L \otimes \mathcal J_\lambda \otimes \cI_\zeta^2$ by $\cO_S$, and complete it to a commutative diagram with exact rows and columns:
\[
\begin{tikzcd}
&  & 0 \arrow{d} & 0 \arrow{d} & \\
0 \arrow{r} & \cO_S \arrow{r}{s_\lambda} \arrow{d}{=} & \tilde E_\lambda \arrow{r} \arrow{d} & L\otimes \mathcal{J}_\lambda\otimes \cI_\zeta^2 \arrow{d} \arrow{r} & 0 \\
0 \arrow{r} & \cO_S  \arrow{r}{s_\lambda'} & \tilde E_\lambda^{**} \arrow{r} \arrow{d} & \tilde E_\lambda^{**}/\cO_S=L\otimes\cI_{Z(s_\lambda')} \arrow{r} \arrow{d} & 0 \\
 &  & \tilde E_\lambda^{**}/\tilde E_\lambda \arrow{d} \arrow{r}{=} & \cI_{Z(s_\lambda')}/\mathcal{J}_\lambda\otimes \cI_\zeta^2 \arrow{r} \arrow{d} & 0 \\
&   & 0 & 0 & 
\end{tikzcd}
\]

Here $s_\lambda$ denotes the section in $H^0(\tilde E_\lambda)$ induced by the extension, and $s_\lambda'$ is the section corresponding to $s_\lambda$ under the inclusion $H^0(\tilde E_\lambda)\subset H^0(\tilde E_\lambda^{**})$. Now observe that $Z(s_\lambda')$ is of local complete intersection at every point; this implies that $Z(s_\lambda')$ has length $\leq 7$ (here we use that any local complete intersection ideal sheaf containing $\cI_\zeta^2$ has colength $\leq 4$). On the other hand, by stability of $E_\lambda^{**}$ we have $7\leq c_2(E_\lambda^{**})=\length Z(s_\lambda')$.

Therefore, we deduce that $Z(s_\lambda')$ has length $7$ and $Z_{cycle}(s_\lambda')=3p+2q+2r$. In particular, $s_\lambda'$ vanishes at $\eta_\lambda\cup\zeta$ for a certain $\eta_\lambda\in\bP(T_pS)$; it follows from the observation in the second paragraph of this proof that $E_\lambda^{**}=E$ (and hence $s_\lambda'\in V^\vee=H^0(E\otimes\cI_{p\cup\zeta})$).

Moreover, for any $\lambda\neq\mu$ and any choice of extensions, the corresponding sections $s_\lambda',s_\mu'$ are not proportional, since they vanish at a different subscheme at $p$. Indeed, any colength 3 complete intersection ideal containing $\mathcal J_\lambda$ must be of the form $(z+\lambda w+\epsilon w^2,w^3)$ for some $\epsilon$. 

Hence varying $s_\lambda$ algebraically with $\lambda$ (and taking closure), we can construct an effective divisor $C\subset \bP(V^\vee)$ of sections whose zero loci are supported in $\{p,q,r\}$.

Now recall that, if $T$ denotes the base locus of the rational map $S\dasharrow \bP(V^\vee)$, then the fiber of $S\setminus T\longrightarrow\bP(V^\vee)$ over $[s]\in\bP(V^\vee)$ is given by $Z(s)\cap (S\setminus T)$ (see \autoref{compdegree}). It follows that the morphism factors through
\[
S\setminus T\longrightarrow \bP(V^\vee)\setminus C,
\]
in particular we have a generically finite morphism from $S$ minus a finite set of points to an affine variety, which is a contradiction.\qedhere
    %
\end{proof}

\section{Genus 13}\label{sec:genus13}
In this section we study a very general polarized $K3$ surface $(S,L)$ of genus $13$ (which in particular implies $\mathrm{Pic}(S)=\mathbb Z\cdot L$). In \cite{mor} it was proven that $\mathrm{irr}_L(S)\le 6$. We are going to improve this bound to $\mathrm{irr}_L(S)\le 4$, by showing that $W^2_4(S,L)$ has a component of positive dimension. 

In this case the moduli space of vector bundles with minimal $c_2=8$ is $\cM=\cM(2,L,6)$, and can be endowed with a primitive polarization $\hat L$ of genus $4$. Note that since $(S,L)$ is very general, $(\cM,\hat L)$ is also very general among the $K3$ surfaces of genus $4$. This enables us to apply a result of Galati-Knutsen (\cite[Proposition 2.2]{GK}), which asserts the existence of (finitely many) curves $C_1,...,C_N\in|\hat{L}|$ with a triple point.

Again by Chen's results \cite{chen} on rational curves, the curves $C_i$ have geometric genus 1 and have the triple point as only singularity. Let us denote by $\nu_i:\tilde C_i\to C_i$ the normalization of $C_i$.
  
\begin{prop}
Let $(S,L)$ be a very general $K3$ surface of genus $13$, and let $C_1,\dots,C_N\in |\hat L|$ be the elliptic curves with a triple point. Then there are nonconstant morphisms
\[
    \tilde C_i\longrightarrow W^2_4(S,L)
\]
($i=1,...,N$) with pairwise disjoint images. In particular, $\mathrm{irr}_L(S)\le 4$ and $\mathrm{dim}(W^2_4(S,L))\ge 1.$
\end{prop}
\begin{proof}
Let $E_i\in \cM$ be the triple point of the elliptic curve $C_i$. The pushforward of any line bundle on $\tilde C_i$ produces a rank 1, torsion-free sheaf on $C_i$, with 3 generators at the singular point $E_i\in C_i$. We thus have a (noncanonical) closed immersion
\[
\tilde C_i\cong J\tilde C_i\hookrightarrow \overline{JC_i}
\]
where $\overline{JC_i}$ denotes the compactified Jacobian of $C_i$ (see \cite[Lemma 3.1]{beauville}), whose image is contained in the locus of torsion-free sheaves with 3 generators at the triple point.

On the other hand, $C_i$ is the unique curve with a triple point at $E_i$ (indeed one can derive a numerical contradiction from $h^1(\cM,\hat{L}\otimes\cI_{E_i})\geq 3$, by considering successive stable extensions). Thus \autoref{coh913} provides an isomorphism
\[
  \{A\in \overline{JC_i}\;|\;h^0(C_i,A\otimes k(E_i))=3\} \cong
\{\xi \in  S^{[4]} \; | \; h^1(S,E_i\otimes \cI_\xi)=3\}.
\]
All in all, we obtain a morphism $\tilde C_i \to W^2_4(S,L)$ (after sending $\xi\mapsto (E_i,H^0(E_i\otimes \cI_\xi))$). This morphism is not constant (otherwise the sections of $H^0(E\otimes \cI_\xi)$ would vanish at infinitely many points). Finally, since the vector bundles $E_i$ are pairwise distinct (recall that there is a unique curve with a triple point at $E_i$), the images of all these morphisms are pairwise disjoint.
\end{proof}

\begin{rem}
If one can prove that $h^0(E\otimes \cI_\zeta)\le 2$ for all $\zeta \in S^{[5]}$, then the map $\tilde C_i\to W^2_4(S,L)$ would be a closed immersion (see also \autoref{singu}).
\end{rem}

\section{Even genera}\label{sec:evengen}
The goal of this section is to prove \autoref{thmB}. This result wants to be a motivation for future work towards a more complete statement. Many proofs are similar (or considerably simpler) to the odd genera case, hence we will provide less details.

As usual, $(S,L)$ will denote a polarized $K3$ surface with $\Pic(S)=\bZ\cdot L$.

\subsection{Genus 8}
Assume that the polarization has genus 8 (i.e. $L^2=14$). Then $\irr_L(S)\leq4$, as proved in \cite{mor}. Indeed, let $E$ denote the stable vector bundle with $v(E)=(2,L,4)$. It satisfies $h^0(E)=6$ and $c_2(E)=5$. Then for every point $p\in S$ we have $h^0(E\otimes\cI_p)=4$, and any $V^\vee\in\Gr(3,H^0(E\otimes\cI_p))$ defines a rational dominant map $\varphi_V:S\dasharrow \bP^2$ of degree $\leq4$.

In other words, we obtain an irreducible component of $W^2_4(S,L)$ as the projectivization of a rank 4 vector bundle on $S$.

In order to determine the equality $\irr_L(S)=4$, we note that according to \autoref{minimalc2}, degree 3 maps correspond to $V^\vee\in\Gr(3,H^0(E))$ such that $\bigcap_{s\in  V^\vee}Z_{cycle}(s)$ is of degree $2$.
The case where the schematic intersection $\bigcap_{s\in  V^\vee}Z(s)$ is of length 2 is discarded by the following:

\begin{lem}\label{vanishingg8}
For every $\xi\in S^{[2]}$, the vanishing $h^1({E}\otimes\cI_\xi)=0$ holds.
\end{lem}
\begin{proof}
We use Bridgeland stability. Both $\cI_\xi$ and $E^\vee[1]$ are $\sigma_{\alpha,-\frac{1}{3}}$-stable objects for every $\alpha>0$, as a consequence of \autoref{minimalrank} (and the fact that there are no rank 3 spherical classes).
In addition, one easily checks that the numerical wall $W(\cI_\xi, E^\vee[1])$ intersects the vertical line $\beta=\frac{-1}{3}$; for the corresponding stability condition $\sigma_0$, both $\cI_\xi$ and $E^\vee[1]$ are $\sigma_0$-stable objects of the same slope, hence $0=\hom(\cI_\xi, E^\vee[1])=h^1(E\otimes\cI_\xi)$.
\end{proof}

On the other hand, if the schematic intersection is a reduced point, then an argument analogous to that of \autoref{vvv} yields a contradiction again. This proves the equality $\irr_L(S)=4$.

\subsection{Genus 10}
Assume now that $L$ is a polarization of genus 10 (i.e. $L^2=18$). Let $E$ denote the stable vector bundle with $v(E)=(2,L,5)$ (which satisfies $h^0(E)=7$, $c_2(E)=6$).

The interesting point for this genus is that $W^2_4(S,L)$ is not irreducible:
\\

\begin{lem}\label{twocomp}
$W^2_4(S,L)$ has a component isomorphic to the Hilbert square $S^{[2]}$, and a component isomorphic to $\mathbb P (E)$.
\end{lem}
\begin{proof}
The component isomorphic to $S^{[2]}$ is elementary to construct: for any $\xi \in S^{[2]}$, the vector space $H^0(E\otimes \cI_\xi)$ has dimension $\geq 3$ and hence induces a map of degree $\leq 4$ by  \autoref{compdegree}. The fact that the map $S^{[2]}\to W^2_4(S,L)$ is an isomorphism onto its image follows from \autoref{singu} together with the remark below the proof of this lemma.

For the second component, pick $[w]\in\bP(E)$, namely $[w] \in \bP(E/E\otimes \mathcal I_p)$ for some $p\in S$. Consider the subsheaf $E_w\subset E$ such that $E_w|_{S\setminus p}\cong E|_{S\setminus p}$, and $E_w=m_p\cdot w+m_p^2E$ in a neighborhood of $p$. It is easy to check that $h^0(E_w)=3$ and that every section in $H^0(E_w)\subset H^0(E)$ vanishes with order $\geq2$ at $p$. Thus by \autoref{compdegree}, $H^0(E_w)$ defines a map of degree $\leq 4$.
\end{proof}

\begin{rem}
Any map of degree 3 is given by $V^\vee\in\Gr(3,H^0(E))$ such that $\bigcap_{s\in  V^\vee}Z_{cycle}(s)$ has degree $3$.

\begin{enumerate}
    \item Bridgeland stability easily gives $h^1(E\otimes \cI_\xi)\le 1$ (hence $h^0(E\otimes \cI_\xi)\le 2$) for any $\xi \in S^{[3]}$. Indeed, a nonzero morphism $\cI_\xi \to E^\vee[1]$ determines the HN filtration of $\cI_\xi$ for suitable stability conditions.
    This discards the case where $\bigcap_{s\in  V^\vee}Z(s)$ has length 3.

    \item If $\bigcap_{s\in  V^\vee}Z(s)$ is a reduced point $p\in S$, then for an appropriate trivialization $E_p\cong e_1\cdot \cO_{S,p}\oplus e_2\cdot \cO_{S,p}$ of $E$ near $p$, the subsheaf $E_1\subset E$ defined by $E_1=m_p\cdot e_1+m_p^3\cdot e_2$ satisfies $h^0(E_1)\geq 3$ (hence $h^1(E_1)\geq 3$). The triple extension
    \[
    0\to\cO_S^{\oplus 3}\to F\to E_1\to 0
    \]
    is stable of rank 5: since $E_1$ has 6 local generators at $p$, it follows that $v(F^{**})=(5,L,k)$ for some $k\geq 2$.

    \noindent If $k\geq 3$, this contradicts the stability of $F^{**}$. For $k=2$, we have that $F$ has exactly 6 local generators at $p$, which implies $h^0(F^{**}\otimes\cI_p)\geq 2$ (hence $h^1(F^{**}\otimes\cI_p)\geq 1$). This easily contradicts the fact that $F^{**}$ is spherical.
\end{enumerate}
The configuration that is hard to exclude is that of a length 2 schematic intersection. Such cases are parametrized by a component of the intersection $\mathbb P(E)\cap S^{[2]}$ of the two loci described in \autoref{twocomp} (it coincides with the intersection along the locus of nonreduced subschemes).
\end{rem}

\subsection{Genus 12 and 14} In this section we consider analogous statements for $K3$ surfaces of genus 12 and 14. In the case of genus 14, the desired component was already described in \cite[Theorem 2.4]{mor}. For genus 12, the following improves the estimate given in \cite{mor}:

\begin{prop}
If $(S,L)$ is a polarized $K3$ surface of genus $12$ with $\Pic(S)=\bZ\cdot L$, then $W^2_4(S,L)$ has a unirational $3$-dimensional component. In particular, $\irr_L(S)\le 4$.
\end{prop}
\begin{proof} 
This can be proved via Bridgeland stability, similarly to other genera. Consider the spherical stable vector bundle $E$ with $v(E)=(2,L,6)$ (so that $c_2(E)=7$ and $h^0(E)=8$).

If $\xi \in S^{[3]}$, any nonzero morphism $\cI_\xi\to E^\vee[1]$ defines a short exact sequence
\[
0\to F\to \cI_\xi\to E^\vee[1]\to 0
\]
(for $F$ the stable bundle with $v(F)=(3,-L,4)$) destabilizing $\cI_\xi$ along a certain wall $W$. This is the HN filtration of $\cI_\xi$ just below $\sigma_0:=W\cap\{\beta=\frac{-2}{5}\}$ (by $\sigma_0$-stability of $F$ and $E^\vee[1]$), hence the morphism  $\cI_\xi\to E^\vee[1]$ is unique (up to constant). It follows that the jump locus $R_3=\{\xi\in S^{[3]}|h^1(E\otimes \cI_\xi)=1\}$ equals $\bP(\Ext^1(E^\vee[1],F))=\bP^3$; for each $\xi$ in this jump locus, $H^0(E\otimes \cI_\xi)$ defines a map of degree $\leq 4$. Hence we get a morphism $\mathbb P^3\to W^2_4(S,L)$. This map is finite onto its image, otherwise we would get $H^0(E\otimes \cI_{\xi_t})$ constant for a $1-$dimensional family of $\xi_t\in R_3$. This would imply that $3=h^0(E\otimes \bigcap \cI_{\xi_t})=h^0(E\otimes \cI_\zeta)=3$ for some $\zeta \in S^{[4]}$ non-curvilinear, thus inducing a map $S\dashrightarrow \bP^2$ of degree $2$, contradiction.
\end{proof}

\section{Final questions and comments}\label{sec:questions}

Throughout this section, we denote by $(S_g,L)$ a very general polarized $K3$ surface of genus $g$.

Even though \autoref{thmA} and \autoref{thmB} are still far from providing a complete description of the number $\irr_L(S_g)$ for arbitrary $g$, they motivate several questions that may deserve some attention. The first natural question is to complete the picture in even genera:

\begin{QU}
 Is $\irr_L(S_g)=4$ for $g=10,12,13,14$? Are there genera $g\geq 15$ for which the equality $\irr_L(S_g)=4$ holds?
\end{QU}

Aside from the independent interest of the Brill-Noether loci, our description seems to indicate that the locus of maps of minimal degree is an essential invariant to understand how $\irr_L(S_g)$ varies with $g$. It is tempting to propose the following:

\begin{QU}\label{quest}
Is it true that either $\irr_L(S_{g+2})=\irr_L(S_g)+1$, or $\irr_L(S_{g+2})=\irr_L(S_g)$ and $\dim W^2_{\irr_L(S_{g+2})}=\dim W^2_{\irr_L(S_{g})}-1$?
\end{QU}

An affirmative answer would imply an affirmative answer to the following question, which is nothing but a weak version of a conjecture in \cite{BDELU}:

\begin{QU}
Does one have $\mathrm{irr}_L(S_g)\to\infty$ as $g\to\infty$?
\end{QU}

Actually, \cite[Conjecture 4.2]{BDELU} predicts something stronger: $\irr(S_g)\to\infty$ as $g\to\infty$. The degree of irrationality of $K3$ surfaces was indeed our original motivation, therefore it is natural to ask:

\begin{QU}
Does the equality $\mathrm{irr}(S_g)=\mathrm{irr}_L(S_g)$ hold for every $g$?
\end{QU}

It is worth mentioning that the analogous question is true for hypersurfaces of large degree in projective space (\cite{BDELU}), and it seems an interesting problem for $(1,d)$-polarized abelian surfaces (for suitable $d$).

All of these questions concern the invariant $\irr_L(S_g)$, but also many questions arise about the Brill-Noether loci, with special focus on $W^2_{\irr_L(S_g)}(S_g,L)$; apart from their dimension, one could ask about other aspects of their geometry (irreducible components, singular locus), Torelli-type questions (e.g. for which $g$ is $(S_g,L)$ determined by $W^2_{\irr_L(S_g)}(S_g,L)$), etc.

In particular, in view of \autoref{singu} it is tempting to draw some parallels about the singularities of the Brill-Noether loci for general $K3$ surfaces with the classical case of curves. As \autoref{singu} suggests, this principle should hold for those components arising as cohomology jump loci. The picture in general seems to be quite complicated to have a comprehensive understanding.

\bibliography{refer}
\bibliographystyle{alphaspecial}

\end{document}